\documentclass{article}

\usepackage{PRIMEarxiv}

\usepackage[utf8]{inputenc} 
\usepackage[T1]{fontenc}    
\usepackage{hyperref}       
\usepackage{url}            
\usepackage{booktabs}       
\usepackage{amsfonts}       
\usepackage{nicefrac}       
\usepackage{microtype}      
\usepackage{lipsum}
\usepackage{graphicx}       
\graphicspath{{media/}}     
\usepackage{algorithmic}
\usepackage{algorithm}%
\usepackage[title]{appendix}%
\usepackage{subfigure}
\usepackage{comment}
\usepackage{wrapfig}
\usepackage{xcolor}
\usepackage{amssymb}
\usepackage{amsmath}
\usepackage{amsthm}
\floatstyle{ruled}
\newfloat{algorithm}{tbp}{loa}
\providecommand{\algorithmname}{Algorithm}
\floatname{algorithm}{\protect\algorithmname}
\theoremstyle{plain}
\newtheorem{thm}{\protect\theoremname}
\theoremstyle{plain}
\newtheorem{assumption}{\protect\assumptionname}
\theoremstyle{plain}
\newtheorem{lem}{\protect\lemmaname}
\theoremstyle{plain}
\newtheorem{rem}{\protect\remarkname}
\theoremstyle{plain}
\newtheorem{dfn}{\protect\definitionname}
\makeatother

\providecommand{\assumptionname}{Assumption}
\providecommand{\lemmaname}{Lemma}
\providecommand{\theoremname}{Theorem}
\providecommand{\remarkname}{Remark}
\providecommand{\definitionname}{Definition}

\newif \iffinal
\finalfalse
\iffinal
\newcommand{\WJmodified}[1]{{#1}}
\newcommand{\WJcomments}[1]{{}}
\newcommand{\oldText}[1]{}
\newcommand{\toAppendix}[1]{}
\else
\newcommand{\WJmodified}[1]{{\color{blue} #1}}
\newcommand{\WJcomments}[1]{{\WJmodified{Wenjie commented: #1}}}
\newcommand{\toAppendix}[1]{#1}
\newcommand{\oldText}[1]{#1}
\fi

\ifpdf
\hypersetup{
  pdftitle={Online Non-convex Optimization with Long-term Non-convex Constraints},
  pdfauthor={Shijie Pan, Jianyu Xu and Wenjie Huang}
}
\fi

\title{Online Non-convex Optimization with Long-term Non-convex Constraints
}

\author{
  Shijie Pan \\
  Department of Electrical and Computer Engineering \\
  Johns Hopkins University \\
  Baltimore, MD 21218, USA\\
  \texttt{span34@jh.edu} \\
   \And
  Jianyu Xu \\
  International Business School Suzhou \\
  Xi'an Jiaotong-Liverpool University \\
  Suzhou, 215123, China\\
  \texttt{Jianyu.Xu@xjtlu.edu.cn} \\
   \And
  Wenjie Huang\thanks{Corresponding author} \\
  Department of Data and Systems Engineering \\
  HKU Musketeers Foundation Institute of Data Science \\
  The University of Hong Kong \\
  Hong Kong SAR, 999077, China\\
  \texttt{huangwj@hku.hk} 
}

\begin{document}
\maketitle

\begin{abstract}
A novel Follow-the-Perturbed-Leader type algorithm is proposed and analyzed for solving general long-term constrained optimization problems in an online manner, where the target and constraint functions are oblivious adversarially generated and not necessarily convex. The algorithm is based on Lagrangian reformulation and innovatively integrates random perturbations and regularizations in primal and dual directions: 1). exponentially distributed random perturbations in the primal direction to handle non-convexity, and 2). strongly concave logarithmic regularizations in the dual space to handle constraint violations. Based on a proposed expected static cumulative regret, and under mild Lipschitz continuity assumption, the algorithm demonstrates the online learnability, achieving the first sublinear cumulative regret complexity 
 for this class of problems.  The proposed algorithm is applied to tackle a long-term (extreme value) constrained river pollutant source identification problem, validate the theoretical results and exhibit superior performance compared to existing methods.
\end{abstract}

\keywords{Online learning,
Long-term constrained optimization, Lagrangian multiplier, Random perturbation, Global minimax point}

\section{Introduction}
%
%
%
%
We study online learning problems where a decision-maker takes decisions over $T$ periods. At each period $t$, the decision $x_{t}\in\mathcal{X\subseteq\mathbb{R}}^{d}$
is chosen before observing a target function $f_{t}$ (not necessarily convex) together with
a set of $I$ constraint functions $c_{it},\,i=1,...,I$ (also not necessarily convex), where the set $\mathcal{X}$ is compact (but not necessarily convex) to guarantee the existence of extreme points. 
Our goal is to find a sequence of decisions $x_t,\,t=1,2,\ldots, T$, whose average target value in the first $T$ periods is as close as possible to that of the best-fixed decision $x^\ast$ in \emph{hindsight}, namely, the decision evaluated through the following offline optimization program that incorporates all observations in the first $T$ periods,
\begin{subequations}
\begin{align}
\min_{x\in\mathcal{X}}\quad & \frac{1}{T}\sum_{t=1}^{T}f_{t}(x)\label{obj}\\
\textrm{s.t.}\quad & \frac{1}{T}\sum_{t=1}^{T}c_{it}(x)\leq b_{i}, & \forall i=1,...,I,\label{con}
\end{align}
\end{subequations}
where each $b_{i}\geq0$ is a certain threshold. This type of online learning problems 
was first explored by Mannor et al. \cite{mannor2009online}, and it has numerous applications ranging from wireless communication \cite{mannor2009online},
GAN network training \cite{razaviyayn2020nonconvex} and repeated auctions \cite{castiglioni2022unifying}. 
Most of the aforementioned studies adopt two primary settings for generating the target and constraint functions: the stochastic setting and the non-oblivious adversarial setting. In the stochastic setting, these functions are drawn from a fixed probability distribution. In contrast, in the non-oblivious adversarial setting, an adaptive adversary generates the functions, possibly in response to the past decisions.
However, in multi-objective online classification~\cite{bernstein2010online}, safe online learning~\cite{amodei2016concrete}, and various sensing applications---such as indoor temperature measurement~\cite{krause2008near}, air quality monitoring~\cite{bai2022novel, qu2017monthly}, and pollution source identification~\cite{huang2023online, wang2018new, yang2016multi, zhang2017pollutant}---the target and constraint functions are typically defined by squared-error losses at observed data points. Observations may be obtained passively—i.e., without intervention
or active querying—or may even arrive in an arbitrary manner.
 In those cases, neither non-oblivious‑adversarial nor stochastic generation setting remain valid. The examples described above illustrate the oblivious adversarial setting \cite{agarwal2019learning,suggala2020online}, where an arbitrary adversary generates the sequence of functions, typically without taking into account the past decisions.
 
Consequently, we consider the following mild assumptions in this paper, where Assumption~\ref{Assu 2} indicates that our setting is the most general among all function‐generation paradigms and thus applies to a much broader class of problems.
\begin{assumption}\label{Assu 1}
    There exists a non-empty compact set where all $I$ constraints in (\ref{con}) hold. Moreover, for any 
 $x\in\mathcal{X}$, each $f_{t}(x)$ and $c_{it}(x),i=1,2,...,I,\,t=1,2,...,T$
satisfy Lipschitz condition with respect to $l_{1}$-norm.
\end{assumption}


\begin{assumption}
\label{Assu 2}
   The target functions (\ref{obj}) and constraint functions (\ref{con}) can be time-varying and
arbitrarily (oblivious adversarially) generated.
\end{assumption}
Any algorithm providing acceptable guarantee against an oblivious adversary inherently maintains its guarantee against non-oblivious adversary or stochastic setting. 
 This establishes that our Assumption \ref{Assu 2} strictly generalizes the non-oblivious adversarial settings considered in prior studies like \cite{castiglioni2022unifying}.

\subsection{Prior Work}
The previous works on online learning for long-term constrained problems are surprisingly sparse. Mannor et al. \cite{mannor2009online} made the first attempt to study online learning problems under sample path constraints for a $2\times 2$ matrix game, where $f_{t}$ and $c_{it}$ are linear functions. They proved that when the opponent (nature) acts non-obliviously, the reward-in-hindsight (i.e., the optimal stationary strategy payoff) is generally unattainable. Nevertheless, they proposed an online algorithm that achieves $O(T^{-1/3})$ average regret relative to the convex hull of the reward-in-hindsight.

Most prior works impose restrictive assumptions. For instance, many studies assume that the constraints are generated i.i.d. from an unknown stochastic model \cite{yu2017online,wei2020online}, or that both the adversarial constraints and target functions satisfy strong structural assumptions, such as smoothness and convexity. In those cases, the generation of the online function sequence is assumed to follow a specific mechanism rather than being oblivious adversarial, or alternatively, a weaker regret metric is adopted \cite{chen2017online,yi2020distributed,cao2018online,sun2017safety}.
 Specifically, Cao and Liu \cite{cao2018online} studied online convex optimization (OCO) with oblivious adversarial generated time-varying convex constraints and convex reward, both of which have uniformly bounded gradients. They derived sublinear regret bounds for both constraint violations and target regrets. The former quantifies the magnitude to which decisions exceed constraint tolerances (LHS of (\ref{con})), and the latter measures the average gap to the optimal constrained target value (\ref{obj}), respectively. The formal definitions of both two regrets will be given in Section \ref{rgf} later.
 Their results extend to bandit feedback settings where the gradient information is unknown but estimable. Chen et al.~\cite{chen2017online} studied a similar setup but adopted a weaker dynamic regret metric, using the myopic optimal solution (considering only the target and constraints at period $t$) as the benchmark, rather than the optimal solution over the entire horizon. Sun et al. \cite{sun2017safety} addressed adversarial contextual bandits with sequential risk constraints. They developed a meta-algorithm leveraging online mirror descent and incorporating expert advice. Yi et al. \cite{yi2020distributed} extended the same OCO to distributed optimization with coupled inequality constraints. Castiglioni et al. \cite{castiglioni2022unifying} unified stochastically
and non-oblivious adversarially generated constraints settings. They proposed an OCO framework for general non-convex functions and arbitrary feasibility sets $\mathcal{X}$. Moreover, they achieved a $\tilde{O}(T^{1/2})$ bound uniformly for the constraint violation regret, target regret and
a $\rho/(1+\rho)$-approximation of hindsight optimum problem (under non-oblivious adversarial constraints). 
Here, $\rho$ is a feasibility parameter related to the existence of strictly feasible solutions. However, their regret bound under stochastically generated setting heavily rely on Azuma-Hoeffding inequality. Under the non-oblivious setting, we should set $\rho=\max_{\xi\in\Xi}\min_{t=1,...,T}\min_{i=1,...,I}-c_{it}(\xi|x_t)$ (constraint function value on strategies mixture $\xi$ is conditioned on real-time decision $x_t$), and it cannot be analyzed by their method under the oblivious adversarial setting. A proximal method of multipliers with quadratic approximations was proposed in \cite{zhang2023regrets} where regrets of the violation of Karush–Kuhn–Tucker (KKT) conditions are analyzed. Under the $L$-smooth, convex compact domain, oblivious adversarial function generation and some other mild conditions, it is shown that this algorithm exhibits $O(T^{-1/8})$ average Lagrangian gradient violation regret, $O(T^{-1/8})$ average constraint violation regret, and $O(T^{-1/4})$ average complementary residual regret if parameters in the algorithm are properly chosen. According to the existing literature, the adversarial cumulative regret bound for online non-convex optimization (ONCO) with the non-obliviously generated rewards and constraints, and under solely the Lipschitz continuity assumption, still remains uninvestigated in the existing literature. 

Meanwhile, unconstrained ONCO has been seen broader progress. Sarkar et~al.\cite{sarkar2025online} studied an online learning problem with a fixed violation budget $B_T$. For an adversarially generated sequence---in which both the $\alpha$-approximately convex objective and constraints are chosen adversarially---of nonnegative, over a convex decision set, they proposed a first-order method that achieves $O(\alpha\sqrt{T})$ on the $\alpha$-static regret.
 Besides, Yang et al. \cite{yang2018optimal} developed a recursive weighting algorithm for an ONCO problem where the solution shall be projected to a non-convex set. Their algorithm retains a $O(d/\sqrt{T})$ average expected regret, where the $d$-dependency in the regret comes from the non-parametric estimation of the non-convex set.
Follow-the-Leader type algorithms \cite{hazan2017efficient,huang2022online} use offline minimization oracle based on gradient descent to attain sublinear regret for a stationary point (not a global minima). Agarwal
et al. \cite{agarwal2019learning} developed a Follow-the-Perturbed-Leader
type (FTPL) algorithm with $O((\mathrm{log}(d)+1)T^{-1/3})$ average expected regret
bound for ONCO problems, though its applicability to long-term constrained problems like (\ref{obj}) warrants further exploration. Unlike standard online algorithms whose regret depends primarily on the diameter of the feasible set, the logarithmic dependence on the dimension $d$ of regret complexity in FTPL arises from the special perturbation scheme. Moreover, the works on unconstrained ONCO typically assume that the target sequence is generated by an oblivious adversary manner, which is the setting adopted in our paper.

\label{section1}
\subsection{Main Contributions of the Paper}
The main contributions of this work are three-unfolded:
\begin{enumerate}
    \item \textbf{Algorithm design}: We propose a novel Follow-the-Perturbed-Leader (FTPL)-type algorithm for online \textit{non-convex} optimization with long-term \textit{non-convex} constraints, where target and constraint functions are generated obliviously and adversarially. The algorithm incorporates two additional terms in each period: a random exponentially distributed linear perturbation \cite{agarwal2019learning} in the primal direction, and a strongly concave logarithmic regularizer in the dual direction. 

\item \textbf{Theoretical guarantee}: We propose an Expected Static Average (ESA)
regret as a comprehensive performance metric for our problem, and derive an $O((\log(d)+1)T^{-1/9})$ average regret bound for the proposed algorithm. To our knowledge, this is the first sublinear regret bound for Problem (\ref{obj})-(\ref{con}) that addresses general non-convexity, incomparable to prior works. Prior studies either focus on stationary-point convergence (first/second-order) \cite{huang2022online, hazan2017efficient, zhang2023regrets} or require restrictive assumptions (e.g., convex domains \cite{zhang2023regrets} and specific constraint-generation mechanisms \cite{castiglioni2022unifying}). From a technical perspective, we extend the online Nash equilibrium framework in \cite{rivera2024online} to a Nash-free setting. We generalize the pure minimization results of \cite{agarwal2019learning} to the primal direction of (\ref{Lag}), while introducing a novel concave dual regularizer for dual complexity estimation. 

\item \textbf{Practical application}: We apply our algorithm to an online river pollutant source identification problem, using streaming sensor data to estimate the released mass, location, and time of the pollutant source upstream. Unlike existing methods \cite{zhang2017pollutant,yang2016multi,wang2018new,huang2023online,bai2022novel}, our approach incorporates long-term (extreme value) constraints to enhance generalization ability. Numerical results demonstrate superior performance of our algorithm, supporting the algorithm’s potential to improve solution reliability in real-world online learning problems under data scarcity.
\end{enumerate}

\subsection{Outline of the Paper}
This paper is organized as follows. In Section 2, we frame the regret formulation to be the performance metric for online learning algorithms. We present the proposed algorithm in Section 3. In the following Section 4, we present the main result on regret complexity bound and then detailed steps of proving the main result. In Section 5, the algorithm is applied to tackle a long-term
(extreme value) constrained river pollutant source identification problem, with superior performance compared to existing methods. Section 6 concludes the paper and points out future research directions.


\section{Regret Formulation}
\label{rgf}
In this section, we propose an Expected Static Average (ESA) 
regret as a comprehensive performance metric for online algorithms solving Problem (\ref{obj})-(\ref{con}), focusing on global optimality rather than local regret measures \cite{zhang2023regrets}.  Existing works \cite{castiglioni2022unifying, zhang2023regrets} measure the performance of online algorithms using either the \textit{constraint violation regret},
\begin{align}
\label{p regret}
\mathfrak{R}_T^{\mathrm{violation}}=\mathbb{E}\left\{\max_{i=1,\,2,\,...,\,I}[\frac{1}{T}\sum_{t=1}^{T}c_{it}(x_{t})-b_i]\right\},
\end{align}
or the \textit{target regret},
\begin{align}
\label{c regret}
    \mathfrak{R}_T^{\mathrm{target}}=\mathbb{E}\left[\frac{1}{T}\sum_{t=1}^{T}f_{t}(x_{t})-\frac{1}{T}\sum_{t=1}^{T}f_{t}(x^\ast)\right].
\end{align}
However, these metrics alone fail to provide a comprehensive evaluation of an algorithm's performance, as they do not simultaneously account for both optimization accuracy and constraint satisfaction. To address this limitation, our study incorporates both constraint violation and target achieving, and defines the period-wise cost for a decision variable $(x,\,y)$ at $t$ as 
\begin{equation*}
L_{t}(x,\,y):=f_{t}(x)+\sum_{i=1}^{I}\gamma_{i}[c_{it}(x)-b_{i}],
\end{equation*}
where $y=(\gamma_{i})_{i=1,...,I}$ are the Lagrangian multipliers penalizing constraint violations. 
To ease the burden of notation, we name the analysis in $x$ as the primal direction and in $y$ as the dual direction. Then, using $L_{t}$, the optimal decision is exactly the solution of the following Lagrangian formulation of Problem (\ref{obj})-(\ref{con}):
\begin{equation}
\label{Lag}
\min_{x\in\mathcal{X}}\max _{y\geq0}\sum_{t=1}^{T}\frac{1}{T}L_{t}(x,\,y)=\min_{x\in\mathcal{X}}\max_{y\geq0}\frac{1}{T}\sum_{t=1}^{T}
\left\{ f_{t}(x)+\sum_{i=1}^{I}\gamma_{i}[c_{it}(x)-b_{i}]\right\},
\end{equation}
namely, by the \emph{Lagrangian Sufficiency Theorem} \cite[Appendix 1, Theorem 5]{courcoubetis2003pricing}, the optimal value of Problem (\ref{Lag}) equals the optimal value of Problem (\ref{obj})-(\ref{con}). 

Our regret analysis (the main result Theorem \ref{thm 1} and its proof shown later) 
can also get rid of checking the Slater's condition for the strict feasibility of the online constrained problem, which is intractable in many real‐world applications, where the underlying feasible set is complicate and may evolve over time. In many existing literature, the Slater's condition is assumed to hold \cite{chen2017online,yi2020distributed,zhang2023regrets}. 
The feasibility of Slater's condition is explicitly verified only in a few restricted scenarios, such as when both the target and constraint functions are generated stochastically~\cite{castiglioni2022unifying}.
Conclusively, by reformulating the constrained optimization into a Lagrangian minimax problem, we can make full use of the rich existing results in unconstrained online games, rather than explicit feasibility checks.

Accounting for randomness in generating solution $\{(x_{t},\,y_{t})\}_{t=1}^{T}$ \cite{rivera2024online,cardoso2019competing}, we propose the ESA 
regret up to period $T$, as,
\begin{equation}
\begin{aligned}
\mathfrak{R}_{T}:=\mathbb{E}\left|\frac{1}{T}\sum_{t=1}^{T}L_{t}(x_{t},\,y_{t})-\min_{x\in\mathcal{X}}\max _{y\geq0}\frac{1}{T}\sum_{t=1}^{T}L_{t}(x,\,y)\right|.\label{oreg}
\end{aligned}
\end{equation}
This regret is adopted by previous works where the offline
oracle is convex in the primal direction and concave in the dual direction, so that global Nash equilibrium \cite{jin2020local} exists. However, given our settings
that both reward and constraint functions are not necessarily convex, the global Nash equilibrium may not exist \cite{jin2020local}. Instead, we consider the \emph{global minimax point} as the hindsight optimal solution of term $\min_{x\in\mathcal{X}}\max _{y\geq0}\sum_{t=1}^{T}L_{t}(x,\,y)$ in (\ref{oreg}). The global minimax point is formally defined as follows.
\begin{dfn} \label{dfn 1} 
\cite{jin2020local}
(Global minimax point) For a function $h$ that satisfies the Lipschitz condition in both primal and dual direction (can be a two-player payoff function) on a compact domain $\mathcal{X}\times\mathcal{Y}$ (not necessarily convex), 
$(x^\ast,\,y^\ast)\in\mathcal{X}\times\mathcal{Y}$ is a global minimax point 
\begin{equation*}
h(x^*,\,y)\leq h(x^*,\,y^*) = \max\limits_{y'\in\mathcal{Y}} h(x^*,\,y')\leq \max\limits_{y'\in\mathcal{Y}} h(x,\,y'),
\end{equation*}
for all $x\in\mathcal{X},\,y\in\mathcal{Y}$.
%
\end{dfn}

\section{The Algorithm and Main Results}

\label{III}
We propose the Follow the Perturbed-Leader (FTPRL) 
algorithm (Algorithm \ref{alg:example}) and present its ESA regret bound. In the inner iteration $m = 1,...,M$, both perturbations and regularizations are imposed. Here i.i.d. random exponentially distributed linear perturbation following $(\mathrm{Exp}(\eta))^d$ with mean $\eta^{-1}$ and variance $\eta^{-2}$ is drawn in the primal direction, following the same idea in \cite{agarwal2019learning}. Strongly concave logarithm function regularizer is imposed to the dual direction. 
For the dual maximization, the original feasible regions is $[0,\,+\infty)$, practical implementation requires bounding the feasible region to $[0,\,y_{\max}]$ where $y_{\max}$ is a sufficiently large positive constant, and then obtain a modified Lagrangian: 
  \begin{equation}
\min_{x\in\mathcal{X}}\max_{y\in[0,y_{\max}]^I}\frac{1}{T}\sum_{t=1}^{T}\left\{ f_{t}(x)+\sum_{i=1}^{I}\gamma_{i}[c_{it}(x)-b_{i}]\right\}. \label{modified Lag}
    \end{equation}
Proposition 2 in the Appendix shows that formulation (\ref{modified Lag}) can serve as a proper approximation for the Lagrangian oracle (\ref{Lag}), when $y_{\max}$ is sufficiently large. So we consider Problem (\ref{modified Lag}) instead of (\ref{Lag}) throughout our regret analysis.
    
In Algorithm \ref{alg:example}, at any period $t$ and iteration $m$, $(x_{tm},\,y_{tm})$ is solved to be the exact global minimax point of the offline oracle $\min_{x\in\mathcal{X}}\max_{y\in[0,y_{\max}]^I} \sum_{n<t}\{\bar{L}_n(x,\,y)-\theta_{t}^{\top}x\}$. 
The solution at period $t$, i.e., $(x_t,\,y_t)$, is then obtained from solving an equation of the sample average of loss function value of $M$ global minimax points. In practice, each offline oracle can be solved by any proper algorithms proposed from the previous studies 
\cite{agarwal2019learning, castiglioni2022unifying, zhang2023regrets}. Specifically, the point $(x_{tm},\,y_{tm})$ can be solved by considering minimizing the following objective,
\[
f(x) = \max_{y\in[0,y_{\max}]^I} \sum_{n<t} \left\{ \bar{L}_n(x,\,y) - \theta_t^{\top} x \right\}.
\]
For any fixed \(x\), the inner maximization with respect to \(y\) is a strongly concave problem on the box domain \([0,\,y_{\max}]^{I}\). It can therefore be solved in polynomial time by \emph{projected gradient ascent}. Once the function value for \(f(x)\) is available, we can obtain an \(\epsilon\)-global minimizer of \(f\) in \(x\) with high probability by using zero-order optimization methods, such as \emph{adapted grid search} \cite{munos2011optimistic, jones2021direct} and \emph{simulated annealing} with a formally justified cooling schedule \cite{hajek1988cooling}. For efficient computational scenarios, one may also employ \emph{genetic algorithms} \cite{mathew2012genetic} or \emph{surrogate-based optimization} techniques \cite{koziel2011surrogate} to efficiently minimize \(f(x)\).

Finally, $L_t(x_t,\,y_t)$ can be solved from the equation-solving step which is equivalent to the optimal solution of following $d + I$ dimentional non-convex square loss minimization problem:
\[
\arg\min_{\substack{x_t \in \mathcal{X},\,y_t \in [0,\,y_{\max}]^I}} 
\left(L_t(x_t,\,y_t) - \frac{1}{M} \sum_{m \leq M} L_t(x_{tm},\,y_{tm})\right)^2,
\]
which 
can also be efficiently solved using the optimization techniques mentioned above. Note that when $\mathcal{X}$ is connected, the the minimal loss can attain zero.
\begin{algorithm}[htbp]
   \caption{Follow the Perturbed-Regularized-Leader (FTPRL)}
   \label{alg:example}
\begin{algorithmic}
\STATE Parameters: $\eta=T^{-2/3}$, $M\geq 1$, and $\lambda > 0$.
\FOR{$t=1,...,T$}
\STATE Construct $\bar{L}_{t}(x,\,y)=L_{t}(x,\,y)+\frac{\lambda I}{t^{1/9}} \sum_{i\leq I}\log(\gamma_i+1).$
\FOR{$m=1,2,...,M$}
\STATE Draw i.i.d random vector $\theta_{t}\sim(\mathrm{Exp}(\eta))^{d}$;
\STATE Solve: $(x_{tm},\,y_{tm})=\arg\,\min_{x\in\mathcal{X}}\max_{y\in[0,y_{\max}]^I} \sum_{n<t}\{\bar{L}_n(x,\,y)-\theta_{t}^{\top}x\}$.
\ENDFOR
\STATE Find $(x_t,\,y_t)$ satisfying equation $L_t(x_t,\,y_t)=\frac{1}{M}\sum_{m\leq M}L_t(x_{tm},\,y_{tm})$.
\ENDFOR
\end{algorithmic}
\end{algorithm}

We distinguish Algorithm \ref{alg:example} from the classical Follow-the-Leader (FTL) algorithm by introducing a primal perturbation term $-\theta_{t}^{\top}x$ and a dual regularizer $\frac{\lambda}{t^{1/9}} \sum_{i\leq I}\log(\gamma_i+1)$. It is well established that online learnability presupposes algorithmic stability: the norm of the update between consecutive iterates must converge to zero. Rivera et al. \cite{rivera2024online} demonstrated that the minimax FTL algorithm 
exhibits stability when the primal-dual function $h(x, y)$ is strongly convex-strongly concave within a compact domain. For the general convex-concave scenario, they introduced two strong convex perturbations separately in primal and dual directions to stabilize the target function. However, these techniques do not extend effectively to the non-convex setting because the Lipschitz continuity of non-convex optimum cannot be guaranteed by a strong convex perturbation. 

Inspired by Agarwal et al. \cite[Lemma 3]{agarwal2019learning}, we demonstrate that the random linear perturbation employed in Algorithm \ref{alg:example} can indeed stabilize the primal direction of $\bar{L}_n(x,y)$, even under non-convex setting. 
Moreover, 
we find that the primal perturbation can also induce stability in the dual direction, provided that $\bar{L}_n(x,y)$ is strongly concave with respect to $y$.

The following Theorem \ref{thm 1} demonstrates the regret complexity for FTPRL, and the online learnability, with the detailed analysis in Section \ref{MR}.
\begin{thm}
\label{thm 1}
Let \( M = \lceil T^{2/9} \rceil \) in Algorithm~\ref{alg:example}. If the decision set \( \mathcal{X} \) is connected, then the ESA regret \( \mathfrak{R}_T\) of FTPRL converges to zero at rate \( O((\mathrm{log}(d)+1)T^{-1/9}) \). 

\end{thm}
\begin{rem}[Implications of Theorem \ref{thm 1}]
\emph{(i) Vanishing rate.} The ESA regret decays polynomially in $T$ and hence
$\mathfrak{R}_T\!\to 0$. The exponent $1/9$ is the current optimal under the analyzing framework FTPRL (see Remark \ref{rem 2} later). To achieve average regret at most $\varepsilon>0$, it suffices that $T=O(((\log d+1)^{9})/\varepsilon^{9})$ and, hiding polylogarithmic factors, $T=\tilde O(\varepsilon^{-9})$. This yields a polynomial sample complexity in $1/\varepsilon$ with exponent $9$. The rate is conservative but guarantees convergence without stronger structural assumptions (e.g., convexity, PL condition), and remains computationally affordable for large-scale problems.

\emph{(ii) Mild dependence on $d$.} The dimension only enters through
$\log d$ (we write $\log d+1$ so the statement is meaningful at $d=1$).
This logarithmic dependence typically arises from the primal perturbation
step in the analysis of FTPRL, and makes the guarantee essentially
dimensional-robust for very large variable spaces.
\end{rem}




\section{Regret Analysis} 
\label{MR}
In this section, we demonstrate the detailed steps and techniques for proving Theorem \ref{thm 1}, through the following subsections.
First, before sketching the intuitive idea of proof, we could first propose another weaker metric as follows.
\[
\mathfrak{R}^{W}_{T}:=\left|\mathbb{E}\left[\frac{1}{T}\sum_{t=1}^{T}L_{t}(x_{t},\,y_{t})-\min _{x\in\mathcal{X}}\max_{y\geq0}\frac{1}{T}\sum_{t=1}^{T}L_{t}(x,\,y)\right]\right|.
\]
We focus on bounding $\mathfrak{R}^{W}_{T}$, as directly analyzing the standard ESA regret $\mathfrak{R}_{T}$ is technically more involved. In particular, our analysis directly bounds the expected regret,
\[
\mathbb{E}\left[\frac{1}{T}\sum_{t=1}^{T}L_{t}(x_{t},\,y_{t}) - \min_{x\in\mathcal{X}}\max_{y\geq 0}\frac{1}{T}\sum_{t=1}^{T}L_{t}(x,\,y)\right],
\]
as a whole, rather than bounding each perturbation parameter $\theta_t$ individually.  
This global bound is made possible by a concentration property of the exponential distribution, namely
$\mathbb{E}[\Vert \theta_t \Vert_{\infty}] \leq \eta^{-1}(\log (d) + 1)$,
which ensures that the magnitude of random perturbations remains controlled.

From the fact of Jensen's inequality, it is known that $\mathfrak{R}^{W}_{T}$ is upper bounded by the ESA regret cumulation $\mathcal{R}_T$ (when the expectation and norm are interchanged). Intuitively, $\mathfrak{R}^{W}_{T}$ could be regarded as a weaker version of ESA regret, which called Weak Expected Static Average (WESA) regret. The WESA regret is not only a strong technical tool in our proof, but also a performance measurement of Algorithm \ref{alg:example} under some extreme conditions as follows.
\newtheorem{pro}{Proposition}
\begin{pro}
\label{pro 2}
For any compact \( \mathcal{X} \), if \( M = 1 \), then the WESA regret \( \mathfrak{R}_T^{W} \) of FTPRL also converges to zero at the same rate \( O((\mathrm{log}(d)+1)T^{-1/9}) \).
\end{pro}
In our analysis, for clarity and tractability, we focus on the cumulative form of the weak ESA regret, namely $T \times \mathfrak{R}^{W}_{T}$, throughout the proof. As outlined above, our first step is to bound the regularized regret expression:
\[
\left|\mathbb{E}\left[\sum_{t=1}^{T}\bar{L}_{t}(x_{t},\,y_{t}) - \min_{x\in\mathcal{X}}\max_{y\geq 0} \sum_{t=1}^{T} \bar{L}_{t}(x,\,y)\right]\right|,
\]
by the cumulative sequential difference in expectation, as established in Lemma~\ref{lem 3}.
Next, we show that both the sequential difference terms and the regularization effects exhibit sublinear growth with respect to the time horizon $T$. In particular, we transform the sequential difference terms from their primal direction into the dual ones (from Lemma~\ref{lem 6} 
to Lemma~\ref{lem 8}) to facilitate bounding. 
After Taylor expansion of $\bar{L}_{t}$ in $x$, we observe a negative quadratic bound on the consecutive dual solution difference, which can be further represented as the difference between Lagrangian values at consecutive primal solutions. Due to the Lipschitz continuity property, this difference is naturally bounded by the consecutive primal solution difference.
Last but not least, with the aid of the positiveness of the dual regularizer, we demonstrate that the accumulated contribution from the regularizer is also sublinear, which allows us to remove the regularizer and conclude the sublinear growth rate of $T \times \mathfrak{R}^{W}_{T}$ (Lemma~\ref{lem 9}). Then, the Proposition~\ref{pro 2} is obtained. Finally, we employ the Sample Average Approximation (SAA) method to derive Theorem~\ref{thm 1} from Proposition~\ref{pro 2}, by repeatedly computing \((x_{tm},\, y_{tm})\) under multiple realizations of primal perturbations in Algorithm \ref{alg:example}.
\subsection{Primal and Dual Lipschitz Continuity}\label{Section 4.1}
Discussion in this section indicates the Lipschitz continuity of $\bar{L}_{n},\,n=1,...,T$ in Algorithm \ref{alg:example}. Referring to the structure of $\bar{L}_n(x,\,y)$, we have:
\[\bar{L}_n(x,\,y)=f_t(x)+\sum_{i=1}^{I}\gamma_{i}[c_{it}(x)-b_{i}]+\sum_{i=1}^{I}\frac{\lambda}{n^{1/9}} \log(\gamma_i+1).\]
Suppose both $f_t(x)$ and $c_{it}(x)$ are $G_0$-Lipschitz continuous (Assumption \ref{Assu 1}). Then by the triangular inequality and the bound of each element in $y$, we have:
\[
\left|\bar{L}_n(x,\,y)-\bar{L}_n(x^{\prime},\,y)\right|\leq (y_{\max}I+1)G_0\Vert x-x^{\prime}\Vert_{1},\]
for all $n=1,2,...,T$, we denote $G_1=(y_{\max}I+1)G_0$. 
And for the dual-direction, we have:
\begin{equation}
\begin{aligned}
&|\bar{L}_n(x,\,y)-\bar{L}_n(x,\,y')|
\\
\leq& \sum_{i=1}^{I}\frac{\lambda}{n^{1/9}} [\log(\gamma_i+1)-\log(\gamma_i'+1)]+\sum_{i=1}^{I}(\gamma_{i}-\gamma_{i}')[c_{it}(x)-b_{i}]\\
\leq& \sum_{i=1}^{I}\frac{\lambda}{n^{1/9}}|\gamma_{i}-\gamma_{i}'|+[c_{it}(x)-b_{i}](\gamma_{i}-\gamma_{i}')\\
\leq& \sum_{i=1}^{I}(\frac{\lambda}{n^{1/9}}+|c_{it}(x)-b_{i}|)|\gamma_{i}-\gamma_{i}'|\leq \Vert y-y'\Vert_1\max_{i=1,2,...,I}\left\{\frac{\lambda}{n^{1/9}}+|c_{it}(x)-b_{i}|\right\}.\nonumber
\end{aligned}
\end{equation}
For the above second inequality, we derive the Lipschitz continuity modulus of $\log(\gamma_i+1),\,\gamma_i\in[0,y_{\max}],\,i=1,2,...,I$, which is obviously upper bounded by 1.
Then, we give the absolute value bound of $[c_{it}(x)-b_{i}](\gamma_{i}-\gamma_{i}')$ for the H\"older's inequality (requirement of non-negativity) to separate two combined terms $|\gamma_{i}-\gamma_{i}'|$ and $\lambda/n^{1/9}+|c_{it}(x)-b_{i}|$. As a consequence, the term, $\max\limits_{i=1,2,...,I}\{\lambda/n^{1/9}+|c_{it}(x)-b_{i}|\}$, is bounded above (by the Lipschitz continuity of each $c_{it}(x),\,i=1,2,...,I,\,t=1,2,...,T$), by a certain value $G_2>0$. So if we denote $G:=\max\{G_1,\,G_2\}$, Assumption \ref{Assu 1} is sufficient to let the Lipschitz condition hold for both primal and dual directions.

\subsection{Bounding a Regularized Form of WESA Regret} \label{Section 4.2}
We fix $M = 1$ in Algorithm \ref{alg:example} and 
define two vector series, $\{y_{t}^{\prime}\}$ and $\{y_{t}^{\prime\prime}\}$, as 
\begin{align*}
y_{t}^{\prime}&:=\arg\max_{y\in[0,\,y_{\max}]^I}\left\{ \sum_{n<t+1}\bar{L}_{n}(x_{t},\,y)-\theta_{t+1}^{\top}x_{t}\right\},\\
y_{t}^{\prime\prime}&:=\arg\max_{y\in[0,\,y_{\max}]^I}\left\{ \sum_{n<t}\bar{L}_{n}(x_{t+1},\,y)-\theta_{t}^{\top}x_{t+1}\right\}.
\end{align*}
Lemma \ref{lem 3} constructs the upper and lower bound of a perturbed form of WESA, which evaluates the performance of the algorithm with respect to $\bar{L}_n$. 
\begin{lem}\label{lem 3} The upper and lower bounds for the perturbed WESA regret w.r.t period $t-1$ are formulated as follows: 
\begin{equation}
    \begin{aligned}
    &\mathbb{E}[\sum_{n<t}\bar{L}_n(x_{n},y_{n})-\min_{x\in\mathcal{X}}\max_{y\in[0,\,y_{\max}]^I}\sum_{n<t}\bar{L}_n(x,y)] 
    \\\leq & G\sum_{n<t}\mathbb{E}[\Vert x_n-x_{n+1}\Vert_1]+\mathbb{E}[\Vert y_n-y_{n+1}\Vert_1]+\mathbb{E}[\Vert y_n^{\prime\prime}-y_{n+1}\Vert_1] 
    \\
    & +2T^{2/3}(\log(d)+1) x_{\max}, \label{upper bound}
\end{aligned}
\end{equation}
    and,
    \begin{equation}
        \begin{aligned}
    &\mathbb{E}[\min_{x\in\mathcal{X}}\max_{y\in[0,\,y_{\max}]^I}\sum_{n<t}\bar{L}_n(x,y)-\sum_{n<t}\bar{L}_n(x_{n},y_{n})]
    \\\leq & G\sum_{n<t}2\mathbb{E}[\Vert x_n-x_{n+1}\Vert_1]+\mathbb{E}[\Vert y_n-y_{n+1}\Vert_1]+\mathbb{E}[\Vert y_{n}^{\prime}-y_{n+1}\Vert_1] 
    \\
    & +2T^{2/3}(\log(d)+1) x_{\max}, \label{lower bound}
\end{aligned}
    \end{equation}
where $x_{\max} := \max_{x\in\mathcal{X}} \|x\|_1$.
\end{lem}
\begin{proof}
We first prove that 
\begin{equation}
  \begin{aligned}
    &\mathbb{E}[\min_{x\in\mathcal{X}}\max_{y\in[0,\,y_{\max}]^I}\{\sum_{n<t}\bar{L}_n(x,y)-\theta^\top_t x\}] \nonumber
    \\
    \geq & \mathbb{E}[\sum_{n<t}\bar{L}_n(x_{n+1},\,y_{n+1})]-G\sum_{n<t}\mathbb{E}[\Vert y_n^{\prime\prime}-y_{n+1}\Vert_1] -T^{2/3}(\log(d)+1) x_{\max}, \label{bridge}
\end{aligned}
\end{equation}
by induction. The first term of LHS and $\bar{L}_1(x_2,y_2)$ of RHS in inequality (\ref{bridge}) can be canceled when $t=2$, and we obtain:
\[G\mathbb{E}[\Vert y_1^{\prime\prime}-y_{2}\Vert_1]+T^{2/3}(\log(d)+1)x_{\max}\geq \mathbb{E}[\theta_{2}^{\top} x_2],\]
which is equivalent to proving the upper bound of RHS, by H{\"o}lder's inequality and the bound of the expected infinite norm on exponential random variable \cite{agarwal2019learning},
\[G\mathbb{E}[\Vert y_1^{\prime\prime}-y_{2}\Vert_1]+T^{2/3}(\log(d)+1)x_{\max}\geq \mathbb{E}[\Vert x_2\Vert_1\Vert \theta_2\Vert_{\infty}].\]
The above inequality naturally holds based on our algorithm settings that,  $\mathbb{E}[\Vert \theta_t\Vert_{\infty}]\leq \eta^{-1}(\log(d)+1)$ for any $t$ and thus $T^{2/3}(\log(d)+1)x_{\max}\geq \mathbb{E}[\Vert x_2\Vert_1\Vert \theta_2\Vert_{\infty}]$. Now we assume inequality (\ref{bridge}) holds some $t$. At $t+1$, we have,
\begin{equation}
\begin{aligned}
&\mathbb{E}[\min\limits_{x\in\mathcal{X}}\max\limits_{y\in[0,\,y_{\max}]^I}\sum_{n<t+1}\bar{L}_n(x,\,y)-\theta^\top_{t+1}x]\\=&\mathbb{E}[\sum_{n<t}\bar{L}_n(x_{t+1},\,y_{t+1})-\theta^\top_{t+1}x_{t+1}+\bar{L}_{t}(x_{t+1},\,y_{t+1})]\\
\geq& \mathbb{E}[\sum_{n<t}\bar{L}_n(x_{t+1},\,y_{t}'')-\theta^\top_{t+1}x_{t+1}]+\mathbb{E}[\bar{L}_{t}(x_{t+1},\,y_{t}'')]
\\=& \mathbb{E}[\sum_{n<t}\bar{L}_n(x_{t},\,y_{t})-\theta^\top_tx_{t}]+\mathbb{E}[\bar{L}_{t}(x_{t+1},\,y_{t}'')],
 \label{bridge 2}
\end{aligned}
\end{equation}
where the last equality holds since  for i.i.d random variable $\theta_t$ in Algorithm \ref{alg:example}, $\mathbb{E}[\theta^\top_tx_{t+1}-\theta^\top_{t+1}x_{t+1}]=0$, and $\mathbb{E}[\theta^\top_tx_{t}-\theta^\top_{t+1}x_{t}]=0$. Combining inequality (\ref{bridge 2}) with inequality (\ref{bridge}) at period $t$, we have
\begin{equation}
\begin{aligned}
&\mathbb{E}[\min\limits_{x\in\mathcal{X}}\max\limits_{y\in[0,\,y_{\max}]^I}\sum_{n<t+1}\bar{L}_n(x,y)-\theta^\top_{t+1}x]
\\
\geq
&\mathbb{E}[\sum_{n<t}\bar{L}_n(x_{n+1},\,y_{n+1})]-G\sum_{n<t}\mathbb{E}[\Vert y_n''-y_{n+1}\Vert_1]\\&-T^{2/3}(\log(d)+1)x_{\max}+\mathbb{E}[\bar{L}_{t}(x_{t+1},\,y_{t}'')]
\\
=&\mathbb{E}[\sum_{n<t+1}\bar{L}_n(x_{n+1},\,y_{n+1})]-G\sum_{n<t}\mathbb{E}[\Vert y_n''-y_{n+1}\Vert_1]+\mathbb{E}[\bar{L}_{t}(x_{t+1},\,y_{t}'')-\bar{L}_{t}(x_{t+1},\,y_{t+1})]
\\
& -T^{2/3}(\log(d)+1)x_{\max}
\\
\geq &\mathbb{E}[\sum_{n<t+1}\bar{L}_n(x_{n+1},\,y_{n+1})]-G\sum_{n<t}\mathbb{E}[\Vert y_n''-y_{n+1}\Vert_1]-G\mathbb{E}[\Vert y_t''-y_{t+1}\Vert_1]
\\
& -T^{2/3}(\log(d)+1)x_{\max}
\\
\geq & \mathbb{E}[\sum_{n<t+1}\bar{L}_n(x_{n+1},\,y_{n+1})]-G\sum_{n<t+1}\mathbb{E}[\Vert y_n''-y_{n+1}\Vert_1]-T^{2/3}(\log(d)+1)x_{\max}, \nonumber
\end{aligned}
\end{equation}
which completes the induction. We also have the following 
\begin{align*}
& \min_{x\in\mathcal{X}}\max_{y\in[0,\,y_{\max}]^I}\sum_{n<t}\bar{L}_n(x,\,y) - \theta_{t}^{\top} x_{t}^{\ast\ast}
\\
\geq & \min_{x\in\mathcal{X}} [\max_{y\in[0,\,y_{\max}]^I}\sum_{n<t}\bar{L}_n(x,\,y) - \theta_{t}^{\top} x]
\\
= & \sum_{n<t}\bar{L}_n(x_{t},\,y_{t})-\theta^\top_{t}x_{t}.\nonumber
\end{align*}
where $x_{t}^{\ast\ast}\in \arg\min_{x\in\mathcal{X}} \max_{y\in[0,\,y_{\max}]^I} \sum_{n<t}\bar{L}_{n}(x,\,y)$.
Thus,
\begin{align*}
&\mathbb{E}[\sum_{n<t+1}\bar{L}_n(x_{n+1},\,y_{n+1})-\min\limits_{x\in\mathcal{X}}\max\limits_{y\in[0,\,y_{\max}]^I}\sum_{n<t+1}\bar{L}_n(x,\,y)]\\
\leq & \mathbb{E}[\sum_{n<t+1}\bar{L}_n(x_{n+1},\,y_{n+1})-\min\limits_{x\in\mathcal{X}}\max\limits_{y\in[0,\,y_{\max}]^I}\sum_{n<t+1}\bar{L}_n(x,y)-\theta^\top_{t+1}x] + \mathbb{E}[\theta_{t+1}^{\top} x_{t+1}^{\ast\ast}]
\\
\leq & G\sum_{n<t+1}\mathbb{E}[\Vert y_{n}''-y_{n+1}\Vert_1]+2T^{2/3}(\log(d)+1)x_{\max},
\end{align*}
where the last inequality holds based on  H{\"o}lder's inequality. From inequality (\ref{bridge}), we have 
\begin{equation}
\begin{aligned}
& \mathbb{E}[\sum_{n<t+1}\bar{L}_n(x_{n},\,y_{n})-\min\limits_{x\in\mathcal{X}}\max\limits_{y\in[0,\,y_{\max}]^I}\sum_{n<t+1}\bar{L}_n(x,\,y)]
\\
\leq & \mathbb{E}[\sum_{n<t+1}\bar{L}_n(x_{n},\,y_{n})]-\mathbb{E}[\sum_{n<t+1}\bar{L}_n(x_{n+1},\,y_{n+1})]
\\
&+G\sum_{n<t+1}\mathbb{E}[\Vert y_n''-y_{n+1}\Vert_1]+2T^{2/3}(\log(d)+1)x_{\max}\\
= & \mathbb{E}[\sum_{n<t+1}\bar{L}_n(x_{n},\,y_{n})]-\mathbb{E}[\sum_{n<t+1}\bar{L}_n(x_{n+1},\,y_{n})]
\\
& +\mathbb{E}[\sum_{n<t+1}\bar{L}_n(x_{n+1},\,y_{n})]-\mathbb{E}[\sum_{n<t+1}\bar{L}_n(x_{n+1},\,y_{n+1})]
\\&
+G\sum_{n<t+1}\mathbb{E}[\Vert y_n''-y_{n+1}\Vert_1]+2T^{2/3}(\log(d)+1)x_{\max}
\\
\leq & G\sum_{n<t+1}\mathbb{E}[\Vert x_n-x_{n+1}\Vert_1]+\mathbb{E}[\Vert y_n-y_{n+1}\Vert_1]+\mathbb{E}[\Vert y_n''-y_{n+1}\Vert_1]+2T^{2/3}(\log(d)+1)x_{\max},\nonumber
\end{aligned}
\end{equation}
where the last inequality holds based on the Lipschitz continuity shown in Section \ref{Section 4.1}, and we complete the derivation of the upper bound (\ref{upper bound}).

Next, we aim to show the lower bound (\ref{lower bound}) and we first prove,
\begin{equation}
\begin{aligned}
&\mathbb{E}[\min\limits_{x\in\mathcal{X}}\max\limits_{y\in[0,\,y_{\max}]^I}\{\sum_{n<t}\bar{L}_n(x,\,y)-\theta^\top_tx\}]
\\
\leq & \mathbb{E}[\sum_{n<t}\bar{L}_n(x_{n+1},\,y_{n+1})]+G\sum_{n<t}\mathbb{E}[\Vert x_n-x_{n+1}\Vert_1]+G\sum_{n<t}\mathbb{E}[\Vert y_{n}'-y_{n+1}\Vert_1]
\\
& +T^{2/3}(\log(d)+1)x_{\max}.
\label{bridge 3}
\end{aligned}
\end{equation}
by induction. First, for the first step ($t=2$), after we eliminate the $\bar{L}(x_{2},\,y_{2})$ on both sides of (\ref{bridge 3}), we obtain,
$-\mathbb{E}[\theta^\top_2x_2]\leq G\mathbb{E}[\Vert x_1-x_{2}\Vert_1]+G\mathbb{E}[\Vert y_{1}'-y_{2}\Vert_1]+T^{2/3}(\log(d)+1)x_{\max}$.
Next, with the same technique for deriving (\ref{bridge 2}), we have,
\begin{equation}
\begin{aligned}
&\mathbb{E}[\min\limits_{x\in\mathcal{X}}\max\limits_{y\in[0,\,y_{\max}]^I}\sum_{n<t+1}\bar{L}_n(x,\,y)-\theta^\top_{t+1}x]\\=&\mathbb{E}[\sum_{n<t}\bar{L}_n(x_{t+1},\,y_{t+1})-\theta^\top_{t+1}x_{t+1}+\bar{L}_{t}(x_{t+1},\,y_{t+1})]
\\\leq& \mathbb{E}[\sum_{n<t}\bar{L}_n(x_{t},\,y_{t}')-\theta^\top_{t+1}x_{t}]+\mathbb{E}[\bar{L}_{t}(x_{t},\,y_{t}')]\\
\leq& \mathbb{E}[\sum_{n<t}\bar{L}_n(x_{t},\,y_{t})-\theta^\top_tx_{t}]+\mathbb{E}[\theta^\top_tx_{t}-\theta^\top_{t+1}x_{t}]+\mathbb{E}[\bar{L}_{t}(x_{t},\,y_{t}')]\\
=& \mathbb{E}[\sum_{n<t}\bar{L}_n(x_{t},\,y_{t})-\theta^\top_tx_{t}]+\mathbb{E}[\bar{L}_{t}(x_{t},\,y_{t}')].\nonumber
\end{aligned}
\end{equation}
Then with the aid of Lipschtiz continuity of $\bar{L}_n(x,\,y)$ (Section \ref{Section 4.1}), we can further the induction for cases $t>2$ as follows,
\begin{equation}
\begin{aligned}
& \mathbb{E}[\min\limits_{x\in\mathcal{X}}\max\limits_{y\in[0,\,y_{\max}]^I}\sum_{n<t+1}\bar{L}_n(x,\,y)-\theta^\top_{t+1}x]
\\
\leq & \mathbb{E}[\sum_{n<t+1}\bar{L}_n(x_{n+1},\,y_{n+1})]+G\sum_{n<t}\mathbb{E}[\Vert x_n-x_{n+1}\Vert_1]+G\sum_{n<t}\mathbb{E}[\Vert y_{n}'-y_{n+1}\Vert_1]+\mathbb{E}[\bar{L}_{t}(x_{t},\,y_{t}')]\\&+T^{2/3}(\log(d)+1)x_{\max}\\
=&
\mathbb{E}[\sum_{n<t+1}\bar{L}_n(x_{n+1},\,y_{n+1})]+G\sum_{n<t}\mathbb{E}[\Vert x_n-x_{n+1}\Vert_1]+T^{2/3}(\log(d)+1)x_{\max}\\&+G\sum_{n<t}\mathbb{E}[\Vert y_{n}'-y_{n+1}\Vert_1]+\mathbb{E}[\bar{L}_{t}(x_{t},\,y_{t}')-\bar{L}_{t}(x_{t+1},\,y_{t}')+\bar{L}_{t}(x_{t+1},\,y_{t}')-\bar{L}_{t}(x_{t+1},\,y_{t+1})]\\
= & \mathbb{E}[\sum_{n<t+1}\bar{L}_n(x_{n+1},\,y_{n+1})]+G\sum_{n<t+1}\mathbb{E}[\Vert x_n-x_{n+1}\Vert_1]+G\sum_{n<t+1}\mathbb{E}[\Vert y_{n}'-y_{n+1}\Vert_1] \nonumber
\\
& +T^{2/3}(\log(d)+1)x_{\max},\nonumber
\end{aligned}
\end{equation}
which completes the induction. We further have 
\begin{align*}
\sum_{n<t}\bar{L}_n(x_{t},\,y_{t})-\theta_{t}^{\top} x_{t}\geq \min_{x\in\mathcal{X}} [\bar{L}_n(x,\,y_{t})] -\theta_{t}^{\top} x_{t} = \min_{x\in\mathcal{X}}\max_{y\in[0,\,y_{\max}]^I}\sum_{n<t}\bar{L}_n(x,\,y) - \theta_{t}^{\top} x_{t},
\end{align*}
and then,
\begin{align*}
& \mathbb{E}[\min_{x\in\mathcal{X}}\max_{y\in[0,\,y_{\max}]^I}\sum_{n<t}\bar{L}_n(x,\,y)]  - \mathbb{E}[\sum_{n<t+1}\bar{L}_n(x_{n+1},\,y_{n+1})]
\\
\leq 
& \mathbb{E}[\min\limits_{x\in\mathcal{X}}\max\limits_{y\in[0,\,y_{\max}]^I}\sum_{n<t+1}\bar{L}_n(x,\,y)-\theta^\top_{t+1}x] - \mathbb{E}[\sum_{n<t+1}\bar{L}_n(x_{n+1},\,y_{n+1})] + \theta_{t}^{\top} x_{t}
\\
\leq & G\sum_{n<t+1}\mathbb{E}[\Vert x_n-x_{n+1}\Vert_1]+G\sum_{n<t+1}\mathbb{E}[\Vert y_{n}'-y_{n+1}\Vert_1]+2T^{2/3}(\log(d)+1)x_{\max}.
\end{align*}
Finally, we have
\begin{equation}
\begin{aligned}
\mathbb{E}&[\min\limits_{x\in\mathcal{X}}\max\limits_{y\in[0,\,y_{\max}]^I}\sum_{n<t+1}\bar{L}_n(x,\,y)-\sum_{n<t+1}\bar{L}_n(x_{n},\,y_{n})]
\\
\leq & \mathbb{E}[\sum_{n<t+1}\bar{L}_n(x_{n+1},\,y_{n+1})-\sum_{n<t+1}\bar{L}_n(x_{n},\,y_{n})]\\&+G\sum_{n<t+1}\mathbb{E}[\Vert x_n-x_{n+1}\Vert_1]+G\sum_{n<t+1}\mathbb{E}[\Vert y_{n}'-y_{n+1}\Vert_1]+2T^{2/3}(\log(d)+1)x_{\max}\\
\leq & G\sum_{n<t+1}2\mathbb{E}[\Vert x_n-x_{n+1}\Vert_1]+\mathbb{E}[\Vert y_n-y_{n+1}\Vert_1]+\mathbb{E}[\Vert y_{n}'-y_{n+1}\Vert_1]+2T^{2/3}(\log(d)+1)x_{\max},\nonumber
\end{aligned}
\end{equation}
which completes the derivation of the lower bound (\ref{lower bound}) and thus completes the proof of Lemma \ref{lem 3}.
\end{proof}

Both the upper and lower bounds in Lemma \ref{lem 3} consist of several sums of consecutive series difference and one sublinear term (comes from the exponentially distributed linear perturbation in the primal direction). Besides, Lemma \ref{lem 3} also shows where the dimension $d$-related term in Theorem \ref{thm 1} comes from. The same technique is also used in the proof of Lemma \ref{lem 7} later on. In the following two subsections, we prove that all of the series sums could be upper-bounded by sublinear complexity  (Lemma \ref{lem 6} and Lemma \ref{lem 8}). 
\subsection{Bounding the Primal Complexity}
The following Lemma \ref{lem 6} shows the sublinear complexity of primal direction $\sum_{n<t}\mathbb{E}[\Vert x_n-x_{n+1}\Vert_1]$. The proof of Lemma \ref{lem 6} relies on existing results in \cite{agarwal2019learning}. We first state the following Lemma \ref{lem 7} originated from \cite{agarwal2019learning}.
\begin{lem}
 \cite[Lemma 7]{agarwal2019learning}
    For any two functions $f_1,f_2:\mathbb{R}^{d}\rightarrow\mathbb{R}$ and vectors $\theta_1,\theta_2\in\mathbb{R}^{d}$, let $x_i(\theta_i)\in\mathop{\arg\min}\{f_i(x)-\theta^\top_ix\},\,i=1,2$. Let $f=f_1-f_2$ and $\theta=\theta_1-\theta_2$, we have that.
    $f(x_1(\theta_1))-f(x_2(\theta_2))\leq\theta^\top(x_1(\theta_1)-x_2(\theta_2)).$ \label{lem 7}
\end{lem}

\begin{lem} \label{lem 6}
    (Primal Complexity) For the $\{x_t\}$ generated by Algorithm \ref{alg:example}, we have:
    \[\sum_{n<t}\mathbb{E}[\Vert x_n-x_{n+1}\Vert_1]\leq O(t\cdot T^{-1/3}).\]
\end{lem}
\begin{proof}
By redefining $\theta:=\theta_t-\theta_{t+1}$, we first prove that $\theta^\top(x_t-x_{t+1})$ has a constant lower bound, which lays the foundation of using Lemma \ref{lem 7}, according to \cite{agarwal2019learning}. We consider the lower bound of $\sum_{n<t+1} \bar{L}_n(x_{t+1},\,y_{t+1})-\sum_{n<t} \bar{L}_n(x_{t+1},\,y_{t}^{\prime\prime}))$ minus the upper bound of $\sum_{n<t+1} \bar{L}_n(x_{t},y_{t}')-\sum_{n<t} \bar{L}_n(x_{t},y_{t}))$ to get the lower bound of $\theta^\top(x_t-x_{t+1})$. Consider the property of the global minimax point, we can prove that,
\begin{equation}
\begin{aligned}
& \theta^\top(x_t-x_{t+1})
\\
\geq & -(\sum_{n<t+1} \bar{L}_n(x_t,\,y_{t}^\prime)-\sum_{n<t} \bar{L}_n(x_t,\,y_{t}))+
(\sum_{n<t+1} \bar{L}_n(x_{t+1},\,y_{t+1})-\sum_{n<t} \bar{L}_n(x_{t+1},\,y_{t}^{\prime\prime}))\\\geq& \bar{L}_t(x_{t+1},\,y_{t}^{\prime\prime})-\bar{L}_t(x_{t},\,y_{t}^{\prime})\geq -2B_L, \label{bound}\nonumber
\end{aligned}
\end{equation}
where $B_L>0$ represents the upper bound of each Lagrangian function value. Hence, if we choose a particular $\theta$, the relationship between $x_{t+1}$ and $x_t$ could be derived. It is easy to obtain that if we exchange the order of the perturbed terms, assign $\theta_1$ to oracle $f_2$ while $\theta_2$ for oracle $f_1$, the exact same relationship will be generated ($f(x_1(\theta_2))-f(x_2(\theta_1))\leq\theta^\top(x_1(\theta_2)-x_2(\theta_1))$). By denoting $\theta_1=\theta_2+3B_L\delta^{-1}\times e_k$, $\delta$ is an scalar and $e_k$ is unit vector with $k$th element as $1$. For any $k \in [d]$, let $x_{k,\min}(\theta) = \min\{x_{t}(\theta)\times e_k,x_{t+1}(\theta)\times e_k\}$, $x_{k,\max}(\theta) = \max\{x_{t}(\theta)\times e_k,x_{t+1}(\theta)\times e_k\}$, then $x_{k,\min}(\theta_1)\geq x_{k,\max}(\theta_2)-\delta$, which is similar to inequality (4) in the proof of Lemma 3 in \cite{agarwal2019learning}. 

Then following the same idea of the proof of Lemma 3 in \cite{agarwal2019learning}, we can prove Lemma \ref{lem 6}. 
\end{proof}
\subsection{Bounding the Dual Complexity}
The complexity of dual direction is shown in Lemma \ref{lem 8}.
\begin{lem} \label{lem 8}
    (Dual Complexity) For the $\{y_t\}$, $\{y_t'\}$ and $\{y_t''\}$ series in Algorithm \ref{alg:example}, we could obtain properties as follows.
    \begin{subequations}
    \begin{align}
    \mathbb{E}[\Vert y_t-y_{t+1}\Vert_1]&\leq O(\sqrt{t^{1/9}T^{-1/3}+(\log(d)+1)t^{-8/9}T^{2/3}}) \label{a}\\
    \mathbb{E}[\Vert y_t-y_t'\Vert_1]&\leq O(t^{-4/9}), \label{b}\\
    \mathbb{E}[\Vert y_{t+1}-y_t''\Vert_1]&\leq O(t^{-4/9}),\label{c}\\
    \mathbb{E}[\Vert y_{t+1}-y_t'\Vert_1]&\leq O(\sqrt{t^{1/9}T^{-1/3}+(\log(d)+1)t^{-8/9}T^{2/3}}),\label{d}\\
    \mathbb{E}[\Vert y_{t}-y_t''\Vert_1]&\leq O(\sqrt{t^{1/9}T^{-1/3}+(\log(d)+1)t^{-8/9}T^{2/3}}).\label{e}
    \end{align}
    \end{subequations}
\end{lem}
\begin{proof}
First, by triangular inequality, bound (\ref{a}) could represent as a combination of bound (\ref{c})  and (\ref{e}) or (\ref{b}) and (\ref{d}). Regarding the function $\sum_{n<t+1} \bar{L}_n(x,y)=\sum_{n<t+1} L_n(x,y)+\sum_{n=1}^{t}\sum_{i=1}^{I}(\lambda I/n^{1/9})\log(\gamma_i+1)$, and referring to the lower bound of the sum of concavity under 1-norm:  $\sum_{n=1}^{t} \lambda I/[n^{1/9}(y_{\max}+1)^2]\geq \lambda_0 I t^{8/9}$, where $\lambda_0=\lambda /(y_{\max}+1)^2$ (This value is the minimum concavity we add in the $t+1$ period), we could obtain an upper bound of each offline oracle by the Taylor expansion with bounded second order term.
\begin{align*}
& \mathbb{E}[\sum_{n<t} \bar{L}_n(x_{t},y)-\theta^{\top}_{t}x_{t}]
\\
\leq& \mathbb{E}[\sum_{n<t} \bar{L}_n(x_{t},y_{t})-\theta^{\top}_{t}x_{t}]-\frac{\lambda_0 I (t-1)^{8/9}}{2}\mathbb{E}[\Vert y-y_{t}\Vert^2_2]
 +\mathbb{E}[(\nabla^{\top}_y\sum_{n<t} \bar{L}_n(x_{t},y_{t})-\theta^{\top}_{t}x_{t})(y-y_{t})]\\
\leq& \mathbb{E}[\sum_{n<t} \bar{L}_n(x_{t},y_{t})-\theta^{\top}_{t}x_{t}]-\frac{\lambda_0 (t-1)^{8/9}}{2}\mathbb{E}[\Vert y-y_{t}\Vert^2_1]
 +\mathbb{E}[(\nabla^{\top}_y\sum_{n<t} \bar{L}_n(x_{t},y_{t})-\theta^{\top}_{t}x_{t})(y-y_{t})],
\end{align*}
where the second inequality is from the norm property that $-\lambda I\Vert y-y_t\Vert_2/(y_{\max}+1)^2\leq-\lambda\Vert y-y_t\Vert_1^2/(y_{\max}+1)^2$. Moreover, the strong concavity parameter, 
\[-\frac{\lambda I(t-1)^{8/9}}{(y_{\max}+1)^2}\geq\nabla^2_{\gamma_i}\left\{\sum_{n<t}\frac{\lambda I}{n^{1/9}}\sum_{i=1}^{I}\log(\gamma_i+1)\right\},\](second order derivative) of the term $\sum_{i=1}^{I}\log(\gamma_i+1)$ ranges within $[-1,-1/(y_{\max}+1)^2]$ for each $i=1,\,2,\,...,\,I$. It means that the Hessian matrix of each offline oracle in Algorithm \ref{alg:example} is not only diagonal but also elementwise upper bounded by negative constant $-\lambda I(t-1)^{8/9}/(y_{\max}+1)^2$ for each $t$. We can choose $\lambda$ as a large number to enlarge the concavity. Due to strong concavity, $\mathbb{E}[(\nabla^{\top}_y\sum_{n<t} \bar{L}_n(x_{t},\,y_{t})-\theta^{\top}_{t}x_{t})(y-y_{t})]\leq 0$ holds for the dual optimum \cite{hazan2016introduction}, so as the global minimax point. The strong concavity with or without primal exponentially distributed linear perturbation is the same because of its independence ($\theta_t^{\top}x_t$ is regarded as a constant), which implies that $\mathbb{E}[(\nabla^{\top}_y\sum_{n<t} \bar{L}_n(x_{t},\,y_{t}))(y-y_{t})]\leq 0$. By rearranging terms, we obtain an inequality as follows:
    \begin{align}
    &\frac{2}{\lambda_0 (t-1)^{8/9}}\mathbb{E}[\sum_{n<t} \bar{L}_n(x_{t},\,y_{t})-\sum_{n<t} \bar{L}_n(x_{t},\,y)]\geq\mathbb{E}[\Vert y-y_{t}\Vert^2_1],\label{inequ 7}
    \end{align}
for all $y\in [0,\,y_{\max}$]. If we consider $\sum_{n<t+1} \bar{L}_n(x,y)$, we could similarly have,
    \begin{align}
    &\frac{2}{\lambda_0 t^{8/9}}\mathbb{E}[\sum_{n<t+1} \bar{L}_n(x_{t+1},y_{t+1})-\sum_{n<t+1} \bar{L}_n(x_{t+1},y)]\geq\mathbb{E}[\Vert y-y_{t+1}\Vert^2_1], \label{inequ 8}
    \end{align}
for all $y\in [0,\,y_{\max}$]. 
 
 To verify the bound (\ref{b}) and (\ref{c}), we can use the property of global minimax point, similar to the proof idea of Lemma \ref{lem 6}, we aim to use a constant to construct upper bounds of $\mathbb{E}[\sum_{n<t} \bar{L}_n(x_{t},y_{t})-\sum_{n<t} \bar{L}_n(x_{t},y_t')]$ in inequality (\ref{inequ 7}),
 \begin{equation}
\begin{aligned}
    & \mathbb{E}[\Vert y_t'-y_{t}\Vert^2_1]
    \\
    \leq & \frac{2}{\lambda_0 (t-1)^{8/9}}\mathbb{E}[\sum_{n<t} \bar{L}_n(x_{t},\,y_{t})-\sum_{n<t} \bar{L}_n(x_{t},\,y_t')]\\
    =&\frac{2}{\lambda_0 (t-1)^{8/9}}\mathbb{E}[\sum_{n<t+1} \bar{L}_n(x_{t},\,y_{t})-\sum_{n<t+1} \bar{L}_n(x_{t},\,y_t')-\bar{L}_{t+1}(x_{t},\,y_{t})+\bar{L}_{t+1}(x_{t},\,y_{t}')]\\
    \leq&\frac{2}{\lambda_0 (t-1)^{8/9}}\mathbb{E}[\sum_{n<t+1} \bar{L}_n(x_{t},\,y_{t}')-\sum_{n<t+1} \bar{L}_n(x_{t},\,y_t')-\bar{L}_{t+1}(x_{t},\,y_{t})+\bar{L}_{t+1}(x_{t},\,y_{t}')]\\
    =&\frac{2}{\lambda_0 (t-1)^{8/9}}\mathbb{E}[\bar{L}_{t+1}(x_{t},\,y_{t}')-\bar{L}_{t+1}(x_{t},\,y_{t})].\nonumber
\end{aligned}
\end{equation}
It is obvious that the expectation part of the last term is upper bound by a constant. By the same reason, the constant upper bound $\mathbb{E}[\sum_{n<t+1} \bar{L}_n(x_{t+1},y_{t+1})-\sum_{n<t+1} \bar{L}_n(x_{t+1},y_t'')]$ in inequality (\ref{inequ 8}) could be derived as follows,
\begin{equation}
\begin{aligned}
    & \mathbb{E}[\Vert y_t''-y_{t+1}\Vert^2_1]
    \\
    \leq&\frac{2}{\lambda_0 t^{8/9}}\mathbb{E}[\sum_{n<t+1} \bar{L}_n(x_{t+1},\,y_{t+1})-\sum_{n<t+1} \bar{L}_n(x_{t+1},\,y_t'')]\\
    =&\frac{2}{\lambda_0 t^{8/9}}\mathbb{E}[\sum_{n<t} \bar{L}_n(x_{t+1},\,y_{t+1})-\sum_{n<t} \bar{L}_n(x_{t+1},\,y_t'')+\bar{L}_{t+1}(x_{t+1},\,y_{t+1})-\bar{L}_{t+1}(x_{t+1},\,y_{t}'')]\\
    \leq&\frac{2}{\lambda_0 t^{8/9}}\mathbb{E}[\bar{L}_{t+1}(x_{t+1},\,y_{t+1})-\bar{L}_{t+1}(x_{t+1},\,y_{t}'')].\nonumber
\end{aligned}
\end{equation}

To verify the bound (\ref{d}), by the property of global minimax point and inequality (\ref{inequ 8}), we have:
\begin{align*}
    & \mathbb{E}[\Vert y_t'-y_{t+1}\Vert^2_1]
    \\
    \leq&\frac{2}{\lambda_0 t^{8/9}}\mathbb{E}[\sum_{n<t+1} \bar{L}_n(x_{t+1},\,y_{t+1})-\sum_{n<t+1} \bar{L}_n(x_{t+1},\,y_t')]\\
    \leq&\frac{2}{\lambda_0 t^{8/9}}\mathbb{E}[\sum_{n<t+1} \bar{L}_n(x_{t},\,y_{t}')-\theta^{\top}_{t+1}x_{t}]-\frac{2}{\lambda_0 t^{8/9}}\mathbb{E}[\sum_{n<t+1} \bar{L}_n(x_{t+1},\,y_t')-\theta^{\top}_{t+1}x_{t+1}]\\
    \leq&\frac{2}{\lambda_0 t^{8/9}}\sum_{n<t+1} G\mathbb{E}[\Vert x_t-x_{t+1}\Vert_1]+\frac{4}{\lambda_0 t^{8/9}}\mathbb{E}[\Vert \theta_{t+1}\Vert_{\infty}]x_{\max}.
\end{align*}
And the term $\mathbb{E}[\sum_{n<t+1}\bar{L}_n(x_{t},y_{t}')-\sum_{n<t+1} \bar{L}_n(x_{t+1},y_t')]$ has an upper bound with respect to $\sum_{n<t+1}\mathbb{E}[\Vert x_t-x_{t+1}\Vert_1]$, due to Lipschitz continuity. Then by the sublinear complexity of the perturbation term by Lemma \ref{lem 6}, we have obtained that $\mathbb{E}[\Vert y_t'-y_{t+1}\Vert^2_1]\leq O(t^{1/9}T^{-1/3}+(\log (d)+1)t^{-8/9}T^{2/3}))$, and thus $\mathbb{E}[\Vert y_t'-y_{t+1}\Vert_1]\leq\sqrt{\mathbb{E}[\Vert y_t'-y_{t+1}\Vert^2_1]}$. The bound (\ref{e}) can be derived in the exact same way as (\ref{d}) as follows,
\begin{equation}
\begin{aligned}
    & \mathbb{E}[\Vert y_t''-y_{t}\Vert^2_1]
    \\
    \leq&\frac{2}{\lambda_0 (t-1)^{8/9}}\mathbb{E}[\sum_{n<t} \bar{L}_n(x_{t},\,y_{t})-\sum_{n<t} \bar{L}_n(x_{t},\,y_t'')]\\
    \leq&\frac{2}{\lambda_0 (t-1)^{8/9}}\mathbb{E}[\sum_{n<t} \bar{L}_n(x_{t+1},\,y_{t}'')-\theta^{\top}_{t}x_{t+1}]-\frac{2}{\lambda_0 (t-1)^{8/9}}\mathbb{E}[\sum_{n<t} \bar{L}_n(x_{t},\,y_t'')-\theta^{\top}_{t}x_{t}]\\
    \leq&\frac{2}{\lambda_0 (t-1)^{8/9}}\sum_{n<t} G\mathbb{E}[\Vert x_t-x_{t+1}\Vert_1]+\frac{4}{\lambda_0 (t-1)^{8/9}}\mathbb{E}[\Vert \theta_{t}\Vert_{\infty}]x_{\max}.\nonumber
\end{aligned}
\end{equation}
Then (\ref{d}) is proved by the Lipschitz continuity and $\mathbb{E}[\Vert y_t''-y_{t}\Vert_1]\leq\sqrt{\mathbb{E}[\Vert y_t''-y_{t}\Vert^2_1]}$. Hence, we finished the proof of Lemma \ref{lem 8}.
\end{proof}
\newtheorem{remark}{Remark}

\subsection{Bridging Regularized WESA, (Original) WESA and ESA Regret Complexity} In this stage, we have shown that:
\begin{equation}
\begin{aligned}
&|\mathbb{E}[\sum_{n<t}\bar{L}_n(x_{n},\,y_{n})-\min_{x\in\mathcal{X}}\max_{y\in[0,\,y_{\max}]^I}\sum_{n<t}\bar{L}_n(x,\,y)]|
\\
\leq & O(\sqrt{t^{19/9}T^{-1/3}+(\log (d)+1)t^{10/9}T^{2/3}}+(\log (d)+1)T^{2/3}+t^{5/9}+tT^{-1/3}),\label{FUZADU}
\end{aligned}
\end{equation}
by combining Lemma \ref{lem 6} and Lemma \ref{lem 8} with Lemma \ref{lem 3}.
However, our ultimate goal is to derive the regret complexity under $L_n(x,y)$, rather than the perturbed form $\bar{L}_n(x,y)$. We find a way to bridge the the perturbed form of WESA and the (original) WESA.
\begin{lem} \label{lem 9}
(Concave-regularizer Relaxation) For both $\{x_t^{\ast\ast}\}$ and $\{y_t^{\ast\ast}\}$ series in which each element is the solution to: 
$\min_{x\in\mathcal{X}}\max_{y\in[0,\,y_{\max}]^I}[\sum_{n<t}L_n(x,\,y)+\frac{\lambda}{n^{1/9}}\sum_{i=1}^{I}\log(\gamma_i+1)]$,
we have:
\begin{align}
&\sum_{n<t}[L_n(x_t^{**},\,y_t^{**})+\frac{\lambda}{n^{1/9}}\sum_{i=1}^{I}\log(\gamma_{it}^{**}+1)]\geq \min_{x\in\mathcal{X}} \max_{y\in[0,\,y_{\max}]^I} \sum_{n<t}L_n(x,\,y),\label{inequ 10}
\end{align}
and,
\begin{align}
&\sum_{n<t}[L_n(x_t^{**},\,y_t^{**})-\frac{\lambda}{n^{1/9}}\sum_{i=1}^{I}\log(\gamma_{it}^{\ddag}+1)]\leq \min_{x\in\mathcal{X}}\max_{y\in[0,\,y_{\max}]^I} \sum_{n<t}L_n(x,y), \label{inequ 12}
\end{align}
where $y_{t}^{\ddag}=\arg\max_{y\in[0,\,y_{\max}]^I}[\sum_{n<t}L_n(x_t^*,\,y)+\frac{\lambda}{n^{1/9}}\sum_{i=1}^{I}\log(\gamma_i+1)]$, and $(x_t^*,\,y_t^*)$ is the solution to: ${\min_{x\in\mathcal{X}}\max_{y\in[0,\,y_{\max}]^I}} \sum_{n<t}L_n(x,\,y)$.
\end{lem} 
\begin{proof}
We define $y_{t}^{\dag}: = {\arg\max_{y\in[0,\,y_{\max}]^I}}\sum_{n<t}L_n(x_t^{**},\,y)$, where $y_{t}^{\dag}:=[\gamma_{1t}^{\dag},\gamma_{2t}^{\dag},...,\gamma_{It}^{\dag}]$. 
 To derive inequality (\ref{inequ 12}), we convert the inequality as follows, 
\begin{equation}
\begin{aligned}
    & \sum_{n<t}[L_n(x_t^{**},\,y_t^{**})+\frac{\lambda}{n^{1/9}}\sum_{i=1}^{I}\log(\gamma_{it}^{**}+1)]
    \\
    \leq & \sum_{n<t}[L_n(x_t^{*},\,y_{t}^{\ddag})+\frac{\lambda}{n^{1/9}}\sum_{i=1}^{I}\log(\gamma_{it}^{\ddag}+1)]\\ \leq & \sum_{n<t}[L_n(x_t^{*},\,y_{t}^*)+\frac{\lambda}{n^{1/9}}\sum_{i=1}^{I}\log(\gamma_{it}^{\ddag}+1)], \label{inequ 11}
\end{aligned}
\end{equation}
and the inequality in (\ref{inequ 11}) holds based on the definition of the global minimax point (like the derivations of (\ref{b}) and (\ref{c}) in Lemma \ref{lem 8}). We prove that inequality (\ref{inequ 12}) by three disjoint cases for elementwise substitution. The basic idea is that our regularization term could be regarded as an extended penalty. Hence, when the 
$i$-th constraint is not satisfied, $\gamma_{it}^{**}=y_{\max}$. And when the constraint is satisfied $\gamma_{it}^{\dag}=0$, then $\gamma_{it}^{**}$ may not be zero because of the existence of the non-negative regularizer. We always have $\gamma_{it}^{**}\geq 
\gamma_{it}^{\dag}$ ($\gamma_{it}^{\dag}$ could either 
 be $0$ (if the corresponding constraint is satisfied) or $y_{\max}$ (if this constraint is violated)) for all $t$ and $i$. Define $\tilde{y}_t$ to be an arbitrary feasible vector in the dual-direction, and $[\tilde{y}_t|\{\tilde{y}_t(i)=c\}]$ implies that the $i$th element of $\tilde{y}_t$ is substituted by a constant $c$. We consider the following cases:
    \subsubsection*{Case 1} $\gamma_{it}^{**}=0$ and $\gamma_{it}^{\dag}=0$, we have \begin{equation}
    \begin{aligned}
    &\sum_{n<t}[L_n(x_t^{**},\,[\tilde{y}_t|\{\tilde{y}_t(i)=\gamma_{it}^{**}\}])+\frac{\lambda}{n^{1/9}}\log(\gamma_{it}^{**}+1)]=\sum_{n<t}L_n(x_t^{**},\,[\tilde{y}_t|\{\tilde{y}_t(i)=\gamma_{it}^{\dag}\}]).\nonumber
    \end{aligned}
    \end{equation}
    \subsubsection*{Case 2} $\gamma_{it}^{**}>0$ and $\gamma_{it}^{\dag}=0$, we have 
    \begin{equation}
\begin{aligned}
&\sum_{n<t}[L_n(x_t^{**},\,[\tilde{y}_t|\{\tilde{y}_t(i)=\gamma_{it}^{**}]\})+\frac{\lambda}{n^{1/9}}\log(\gamma_{it}^{**}+1)]
\\
\geq & \sum_{n<t}L_n(x_t^{**},\,0)=\sum_{n<t}L_n(x_t^{**},\,[\tilde{y}_t|\{\tilde{y}_t(i)=\gamma_{it}^{\dag}]\}).\nonumber 
\end{aligned}
\end{equation}
    The inequality of Case 2 holds because $\gamma_{it}^{**}$ maximizes the function $\sum_{n<t}\bar{L}_n(x_t^{**},y)$, hence it is larger than or equal to the situation of $\gamma=0$. When $\gamma=0$, the regularization term $\lambda\log(\gamma+1)/n^{1/9}$ is equal to zero and thus disappeared. That is, $\sum_{n<t}\bar{L}_n(x_t^{**},0)=\sum_{n<t}L_n(x_t^{**},0)$.
\subsubsection*{Case 3} $\gamma_{it}^{**}=\gamma_{it}^{\dag}=y_{\max}$, we have \begin{align*}
&\sum_{n<t}[L_n(x_t^{**},\,[\tilde{y}_t|\{\tilde{y}_t(i)=y_{\max}\}])+\frac{\lambda}{n^{1/9}}\log(y_{\max}+1)]
\geq \sum_{n<t}L_n(x_t^{**},\,[\tilde{y}_t|\{\tilde{y}_t(i)=y_{\max}\}]).
    \end{align*}
    When $\gamma_{it}^{\dag}=y_{\max}$, it means that $\sum_{n<t}[c_{it}(x)-b_{i}]\geq 0$, so that the dual subproblem of offline oracle corresponding to any $i$, is monotonically increasing when $\gamma_i$ becomes larger. Then $\gamma_{it}^{**}=y_{\max}$. By accumulating the above LHS with all $i=1,...,I$, for all constraints, 
    \begin{equation}
        \begin{aligned}
              &\sum_{n<t}[L_n(x_t^{**},\,y_t^{**})+\frac{\lambda}{n^{1/9}}\sum_{i=1}^I\log(\gamma_{it}^{**}+1)]\\&\geq\sum_{n<t}[L_n(x_t^{**},\,[y_t^{**}|\{y_t^{**}(1)=\gamma_{1t}^{\dag}\}])
            +\frac{\lambda}{n^{1/9}}\sum_{i=2}^I\log(\gamma_{it}^{**}+1)]\\&\geq \sum_{n<t}[L_n(x_t^{**},\,[y_t^{**}|\{y_t^{**}(1)=\gamma_{1t}^{\dag}\}\cap\{y_t^{**}(2)=\gamma_{2t}^{\dag}\}]+\frac{\lambda}{n^{1/9}}\sum_{i=3}^I\log(\gamma_{it}^{**}+1)])
              \\
              & ...
              \\&\geq\sum_{n<t}L_n(x_t^{**},\,[y_t^{**}|\{y_t^{**}(1)=\gamma_{1t}^{\dag}\}\cap...\cap \{y_t^{**}(I)=\gamma_{It}^{\dag}\}])\geq\sum_{n<t}L_n(x_t^{**},\,y_{t}^{\dag})\geq\sum_{n<t}L_n(x_t^{*},\,y_{t}^*),\nonumber
        \end{aligned}
    \end{equation}thus we obtain inequality (\ref{inequ 10}). In the proof, we use $y_t^{**}$ initially, and make elementwise substitution by $\gamma_{it}^{\dag}$ from $i=1$ to $I$ and finally obtain the term with $y_{t}^{\dag}$. Here $[\tilde{y}_t|\{\tilde{y}_t(i)=c_1\}\cap \{\tilde{y}_t(j)=c_2\}]$ implies that the $i$th and $j$th elements of $\tilde{y}_t$ are substituted by constant $c_1$ and $c_2$, respectively, and so on.
    Hence, the proof of Lemma \ref{lem 9} is completed. 
    \end{proof}

By combining the results of bound (\ref{FUZADU}) and Lemma \ref{lem 9}, the global minimax point of $\sum_{n<t}\bar{L}_n(x,\,y)$ can be transformed to the one of  $\sum_{n<t}L_n(x,\,y)$ by adding sublinear complexity terms. 
\begin{equation}
\begin{aligned}
&\mathbb{E}[\sum_{n<t}L_n(x_{n},\,y_{n})-\min_{x\in\mathcal{X}}\max_{y\in[0,\,y_{\max}]^I}\sum_{n<t}L_n(x,\,y)]
\\
& -\frac{\lambda}{n^{1/9}}\log(y_{t}^{\ddag}+1)+\frac{\lambda}{n^{1/9}}(\log(y_n+1)+\log(y_t^{**}+1))\\\leq & \mathbb{E}[\sum_{n<t}\bar{L}_n(x_{n},\,y_{n})-\min_{x\in\mathcal{X}}\max_{y\in[0,\,y_{\max}]^I}\sum_{n<t}\bar{L}_n(x,\,y)]
\\
\leq & O(\sqrt{t^{19/9}T^{-1/3}+(\log (d)+1)t^{10/9}T^{2/3}}+t^{5/9}+tT^{-1/3}+(\log (d)+1)T^{2/3}).\nonumber
\end{aligned}
\end{equation}
Hence, by rearranging terms, we have:
\begin{equation}
\begin{aligned}
&\mathbb{E}[\sum_{n<t}L_n(x_{n},\,y_{n})-\min_{x\in\mathcal{X}}\max_{y\in[0,\,y_{\max}]^I}\sum_{n<t}L_n(x,\,y)]
\\
\leq & O(\sqrt{t^{19/9}T^{-1/3}+(\log (d)+1)t^{10/9}T^{2/3}}+tT^{-1/3}+t^{5/9}+(\log (d)+1)T^{2/3}+t^{8/9}).\nonumber
\end{aligned}
\end{equation}
Based on the bound (\ref{FUZADU}), we can also derive the same complexity bound for: \[\mathbb{E}[\min_{x\in\mathcal{X}}\max_{y\in[0,\,y_{\max}]^I}\sum_{n<t}L_n(x,\,y)-\sum_{n<t}L_n(x_{n},\,y_{n})],\] 
and obtain,
\begin{align*}
T\times\mathfrak{R}_{T}^{W} = & \left|\,\mathbb{E}\left[\min_{x\in\mathcal{X}} \max_{y\in[0,\,y_{\max}]^I} \sum_{n<t} L_n(x,\,y)
- \sum_{n<t} L_n(x_n,\,y_n)\right]\,\right| \\
\leq &\; O\Bigl( (\log (d) + 1)\bigl( 
\sqrt{t^{19/9}T^{-1/3} + t^{10/9}T^{2/3}} 
+ tT^{-1/3} + t^{5/9} + T^{2/3} + t^{8/9} 
\bigr) \Bigr).
\end{align*}
By assigning $t=T$, the RHS of the above inequality becomes $ O((\mathrm{log}(d)+1)T^{8/9})$. Both the perturbation and regularization terms have the order lower or equal to $O((\mathrm{log}(d)+1)T^{8/9})$ and will not affect the overall complexity.
Here we derive the WESA regret complexity by dividing $T$ over the growth rate as follows. 

Then, we can utilize the sample average approximation (SAA) method (set $M>1$ in Algorithm \ref{alg:example}) to estimate each expected loss function value $\mathbb{E}[L_t(x_t,\,y_t)]$ to link ESA regret complexity $\mathfrak{R}_{T}$ (Theorem \ref{thm 1}) with the WESA one (Proposition~\ref{pro 2}). To analyze the complexity brought by SAA,  
    the idea is that the average online oracle solution could be considered as an unbiased estimator to $\mathbb{E}[L_n(x_n,\,y_n)]$ for the $M=1$ case. With the aid of triangular inequality, we obtain,
   \begin{align*}
    T\times\mathfrak{R}_{T} 
    \leq & \mathbb{E}[\sum_{n<T}|L_n(x_n,\,y_n)-\frac{1}{M}\sum_{m\leq M}L_n(x_{nm},\,y_{nm})|]
    \\
    & +\mathbb{E}[\sum_{n<T}|\frac{1}{M}\sum_{m\leq M}L_n(x_{nm},\,y_{nm})-\mathbb{E}[L_n(x_n,\,y_n)]|]+T\times\mathfrak{R}_{T}^{W}.
   \end{align*}
    Suppose there is no error comes from the equation solution in Algorithm \ref{alg:example}, the term $\mathbb{E}[\sum_{n<T}|\frac{1}{M}\sum_{m\leq M}L_n(x_{nm},\,y_{nm})-\mathbb{E}[L_n(x_n,\,y_n)]|]$ then has complexity $O(M^{-1/2}T)$ \cite{roy2019multi}. Summarizing the two parts, we obtain the complexity bound on ESA regret as,
    $\mathfrak{R}_{T} \leq O(M^{-1/2}+(\mathrm{log}(d)+1)T^{-1/9})$.
    And if we set $M=\lceil{T^{2/9}}\rceil$, 
    the proof of Theorem \ref{thm 1} is completed. 

\begin{rem}
(Why $O((\mathrm{log}(d)+1)T^{-1/9})$ is the best order to achieve ?) Suppose that the offline oracle invoked at each iteration satisfies an $O((\mathrm{log}(d)+1)T^k)$-concavity condition at step $T$.  
Then, according to Lemma~\ref{lem 8}, the overall complexity of the WESA regret is bounded by,
\[
O((\mathrm{log}(d)+1)(T^{-1/3} + T^{k-1})) + O\left((\mathrm{log}(d)+1)\sqrt{T^{2/3} \cdot T^{-k}}\right).
\]
Following the balancing principle introduced in~\cite{baby2024online}, the minimal regret bound is achieved by choosing $k$ such that the two dominant terms are of the same order.  
That is, we solve for the balancing condition $T^{k-1} = T^{1/3 - k/2}$.
Hence, we have $k=8/9$, which guides us setting the modulus $\lambda/n^{1/9}$ in dual strongly concave logarithm regularizer.
\label{rem 2}
\end{rem}

\section{Numerical Experiments and Applications}
We apply the proposed FTPRL algorithm (Algorithm \ref{alg:example}) to a river pollutant source identification problem \cite{huang2023online} in online manner, by streamingly observed concentration levels in each period. According to existing literature, we use the advection-dispersion equation (ADE) model to simulate the pollutant diffusion process in the river \cite{wang2018new}, 
\begin{equation}
\begin{aligned}
     C(l,\,t|s_0,\,l_0,\,t_0)
  :=\frac{s_0\exp(-k(t-t_0))}{A\sqrt{4\pi D(t-t_0)}}\exp\left\{-\frac{(l-l_0-v(t-t_0))^2}{4D(t-t_0)}\right\}, \label{ADE}
\end{aligned}
\end{equation}
where $s_0,\,l_0,\,t_0$ represent the pollutant source information: the mass of pollutant, the location of source, and the released time. ADE model basically estimates the pollutant concentration level at downstream location $l$ and time $t$, given the known source information $s_0,\,l_0$, and $t_0$. 

For the parameters in the ADE model, we choose the parameters: $D=2430m^2/min$, dispersion coefficient; $k=0min^{-1}$, decay coefficient; $A=60m^2$, area perpendicular to the river flow; $v=80m/min$, the velocity of the river flow, according to the case of Rhodamine WT dye concentration data from a travel time study on the Truckee River between Glenshire Drive near Truckee, Calif., and Mogul, Nev., June 29, 2006 \cite{crompton2008traveltime}. The same parameter settings are also used in \cite{huang2023online, wang2018new}.

Our goal is to use the ADE model to estimate the source information $s_0,\,l_0,\,t_0$, given the streaming concentration data $c_n^o$ at $n$th sampling collected by a sensor at location $l^n$ and collection time $t^n$. The performance of regression is measured by the gap between the simulation value from the ADE model by fixing the source information to be $(s,\,l,\,t)$ and real data $c_n^o$. For simplicity, we use $c_n^p(s,\,l,\,t)$ to represent $C(l^n,\,t^n|s,\,l,\,t)$. The source identification can be formulated into a long-term constrained optimization problem with square error minimization objective and extreme value constraint, over the whole $T$ periods.
The constraint here controls the value of extremely large square loss, which improves the generalization ability of the identification results,
\begin{align}
\min_{x\in\mathcal{X}}\quad & \sum_{n=1}^T(c_n^p(s,\,l,\,t) - c_n^o)^2\nonumber
\\
\quad \textrm{s.t.} \,\,\, &\frac{1}{T}\sum_{n=1}^T \exp((c_n^p(s,\,l,\,t) - c_n^o)^2)\leq b, \nonumber
\end{align}
where $x = (s,\,l,\,t)$ and $\mathcal{X}$ is a box set, $\mathcal{X}\::=\;[10^6,1.5\times10^6]\times[-25000,-20000]\times[-250,-200]$. 
Extremely large error points are penalized via an exponential loss penalty function, and controlled by a properly chosen threshold $b=1.3$.  The experiments
mainly contain three parts: the regret analysis, the out-of-sample performance analysis, and the identification accuracy
analysis. In terms of data generation, the value $l^n$ are chosen from the $30$ identical sections ($30$ different sensors) within the interval $l^n\in[-14216.3,\,22009]$ from Rhodamine WT dye case and $c_n^o$ is generated through ADE model with an underlying ``true'' source information $[s,\,l,\,t] = [1300000,-22106,-215]$ with normally distributed $\mathcal{N}(0,\,0.5)$ random errors, similar to the settings in \cite{huang2023online}. To validate the theoretical regret bound of our algorithm, we plot the average cumulative regret $\mathfrak{R}_{n}/n$ shown in Figure \ref{fig.1}, which shows an obvious descending trend by choosing $\eta = 500^{2/3},\,\lambda=100$, $y_{\max}=100$ and $M=1$. 

When implementing Algorithm \ref{alg:example}, we choose Genetic Algorithm (GA) to solve the offline minimax oracle in the experiments because it gives the best performance among all heuristics. 
\subsection{Regret Analysis}
\begin{figure}[ht]
    \centering
    \subfigure{
        \includegraphics[width=2.9in]{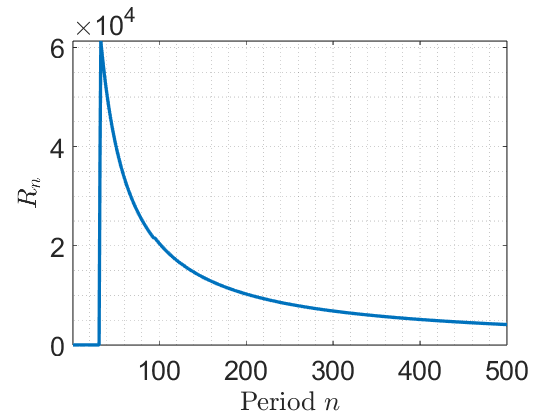}
        \label{label_for_cross_ref_1}
    }
    \subfigure{
	\includegraphics[width=2.9in]{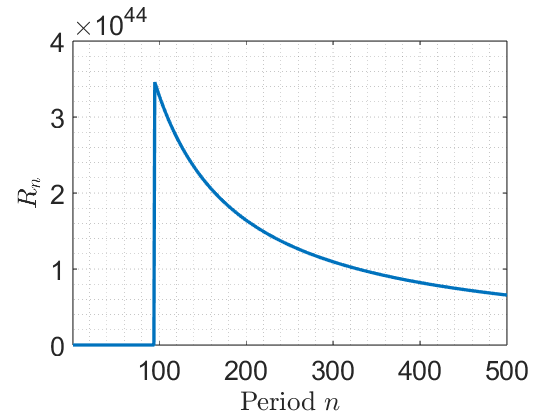}
        \label{label_for_cross_ref_2}
    }
    \caption{Average cumulative regret of Algorithm \ref{alg:example} (left) and Baseline (right)}
    \label{fig.1}
\end{figure}

\begin{figure}[ht]
    \centering
     
    \subfigure{
    	\includegraphics[width=2.9in]{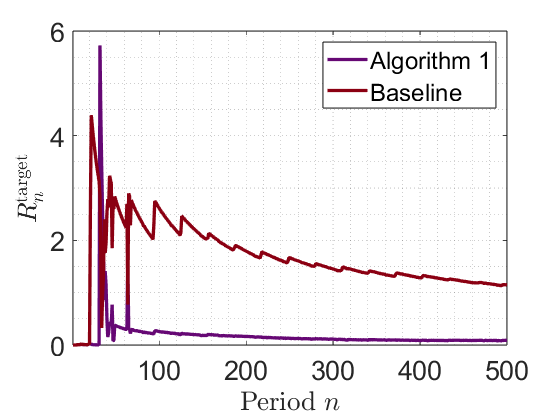}
        \label{label_for_cross_ref_3}
    }
    \subfigure{
	\includegraphics[width=2.9in]{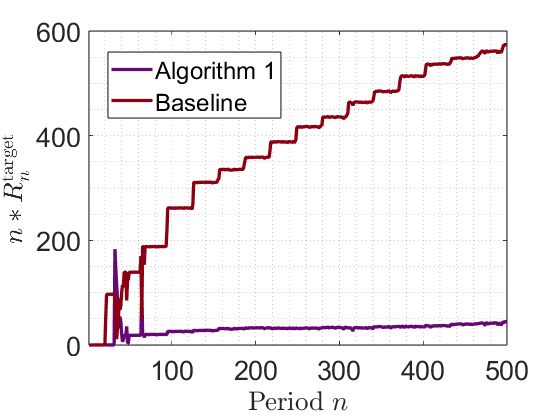}
        \label{label_for_cross_ref_4}
    }
    \caption{Average (left) and Cumulative (right) target regret (\ref{p regret}) of Algorithm \ref{alg:example} (purple) and Baseline (red) 
    }
    \label{fig.3}
\end{figure}

\begin{figure}[ht]
    \centering
     
    \subfigure{
    	\includegraphics[width=2.9in]{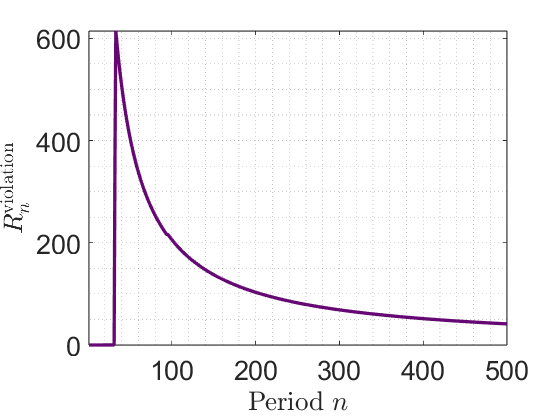}
        \label{label_for_cross_ref_31}
    }
    \subfigure{
	\includegraphics[width=2.9in]{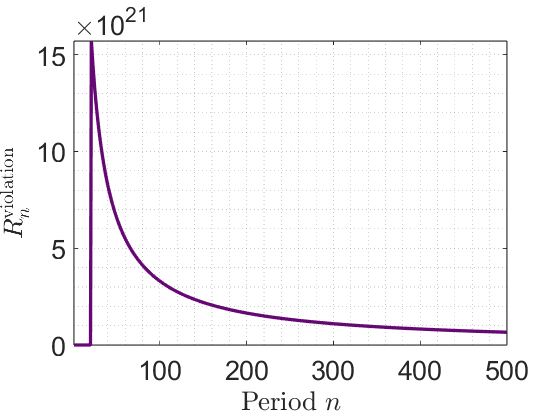}
        \label{label_for_cross_ref_41}
    }
    \caption{Cumulative constraint violation regret (\ref{c regret}) of Algorithm \ref{alg:example} (left) and Baseline (right) 
    }
    \label{fig.5}
\end{figure}

In the first experiment, we implement our algorithm (compared with the baseline \cite[Algorithm 1]{zhang2023regrets}) and plot the first $500$ periods for the average cumulative regret in Figure \ref{fig.1}.
Note that the baseline algorithm in \cite{zhang2023regrets} does not impose an upper bound on each penalty factor $y_i$, so the average cumulative regret is compatible to analyze the baseline algorithm. 
Besides, the target and constraint violation regret defined in \cite[Section 2.2]{castiglioni2022unifying} (also see them in the Introduction section) 
are shown in Figure \ref{fig.3} and \ref{fig.5}, respectively. Both the graphs of average target regret and average constraints violation regret rapidly diminish to zero. Our experiment results support our theoretical results of Theorem \ref{thm 1} and show our regret complexity advantage over the results in \cite{zhang2023regrets}.
  
\subsection{Out-of-Sample Performance Analysis}
\begin{figure}[ht]
\centering  
\includegraphics[scale=0.7]{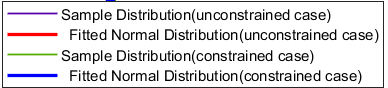}
\\
\includegraphics[scale=0.4]{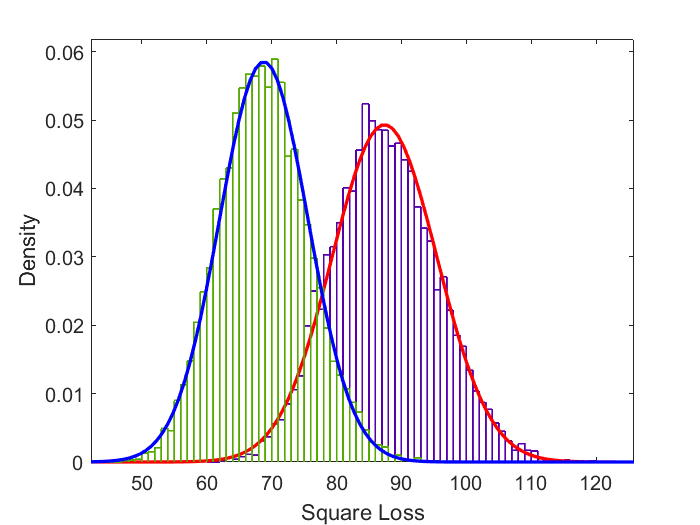}
\caption{Out-of-Sample Performance for constrained ($b=1.3$) and unconstrained problem, 
corresponds to the cross-section period 400 
}
\label{Appendix: fig}
\end{figure}

In the second experiment, we illustrate that the out-of-sample test performance of our solutions is superior to those from the unconstrained case \cite[Algorithm 1]{agarwal2019learning} and the constrained case \cite[Algorithm 1]{zhang2023regrets}. We take five different cross-section periods, to evaluate the mean and variance of square loss generated by online identification solutions of all three algorithms (see Table \ref{Table 1}). Specifically, we sum of all observed losses until all different cross-
section periods, e.g. for the period $400$, we sum up the pointwise square loss from period $1$ to $400$ for evaluation. It demonstrates the merit of extreme value constraint that our solutions yield lower mean loss and variance across all the cross-section cross-
sections. As an example, the one with cross-
sections 400 is explicitly shown by square loss distributions in Figure \ref{Appendix: fig}.

\begin{table}[ht] 
    \caption{Out-of-Sample Performance Comparison} \label{Table 1}
    \centering
    \begin{tabular}{lllll}
    \hline
      \textbf{Period} & & \textbf{Algorithm \ref{alg:example}} & \textbf{\cite[Algorithm 1]{agarwal2019learning}} & \textbf{\cite[Algorithm 1]{zhang2023regrets}}\\
      \hline
      $500$ &$\mu$& $91.2395$ & $106.9429$ & $556.6771$\\
        &$\sigma^2$&  $63.5185$ & $79.9949$ & $513.3791$\\
        \hline
      $450$ &$\mu$& $78.8883$ & $105.2839$ & $517.274$\\
      &$\sigma^2$& $55.2465$ & $78.7825$ & $484.8949$\\
      \hline
      $400$ &$\mu$ & $68.7711$ & $87.5011$ & $475.5414$\\
      &$\sigma^2$& $46.8819$ & $54.3950$ & $440.4379$\\
      \hline
      $350$ &$\mu$& $59.8691$ & $74.8828$ & $195.9447$\\
      &$\sigma^2$& $39.9766$ & $54.5316$ & $176.8275$\\
      \hline
      $300$ &$\mu$ & $51.7917$ & $64.2299$ & $395.2411$\\
      &$\sigma^2$ & $35.9434$ & $49.0830$ & $378.6152$\\
      \hline
    \end{tabular}
  \label{table 1}
\end{table}
 The possible reason of the lower mean, could be explained by a weight perspective. Given a certain $b$, the solution feasible to our offline target: $\min _{x\in\mathcal{X}}\max_{y\in[0,\,y_{\max}]^I}\sum_{t=1}^{T}L_{t}(x,\,y)$, could be infeasible to offline oracle because of perturbation and regularization terms. Hence, all $y_t$ reach their upper bound $y_{\max}$, which contribute to a large value in penalty terms $y\exp\left[(c_{n}^{p}(s,\,l,\,t)-c_{n}^{o})^{2}\right]-Tb$. So the influence of perturbation and regularization terms (particularly the exponentially distributed linear perturbation in the primal direction) will not much influence the offline oracle optimization problem. 
On the other hand, the exponentially distributed linear perturbation in unconstrained variant \cite[Algorithm 1]{agarwal2019learning} will influence the solution by a great extent.
\subsection{Identification Accuracy Analysis}
The relative error of identification solutions on each cross-section is shown in Table \ref{table 2}. We find that the identification result becomes more accurate as the learning process evolves. In addition, we compare the results of our Algorithm \ref{alg:example}, with \cite[Algorithm 1]{agarwal2019learning} and \cite[Algorithm 1]{zhang2023regrets}, which shows our apparent superior performance on $l$ and $t$ directions. Note that the identification result by \cite[Algorithm 1]{zhang2023regrets} seems to converge to the boundary of feasible region $[1000000, -20000, -200]$. Thus it fails in this problem, so we do not show its result in Table \ref{table 2}. 
\begin{table}[H]
    \caption{Identification Result} 
    \centering
    \label{table 2}
    \begin{tabular}{lll}
    \hline
      \textbf{Period} & \textbf{Algorithm \ref{alg:example}} & \textbf{\cite[Algorithm 1]{agarwal2019learning}}\\
      \hline 
      $10$ & $[1488100,-20270,-201]$ & $[1495004, -20000, -200]$\\
      $30$ & $[1489400,-20332,-205]$ & $[1493296, -20000, -200]$\\
      $450$ & $[1456700,-22387,-219]$ & $[1478600,-21118,-202]$\\
      $500$ & $[1470600,-22397,-219]$ & $[1461100,-20789,-200]$\\
      \hline
      \textbf{Period} & \textbf{Relative Error} & \textbf{Relative Error (Baseline)} \\
      \hline 
      $10$ & $[14.47\%,8.31\%,6.67\%]$ & $[15\%,9.53\%,6.98\%]$\\
      $30$ & $[14.57\%,8.02\%,4.54\%]$ & $[14.87\%,9.53\%,6.98\%]$\\
      $450$ & $[12.05\%,1.27\%,1.8\%]$ & $[13.74\%,4.47\%,6.05\%]$\\
      $500$ & $[13.12\%,1.32\%,1.83\%]$ & $[12.39\%,5.96\%,6.98\%]$\\
      \hline
    \end{tabular}
\end{table}

\section{Conclusion}
We develop and analyze a novel algorithm (FTPRL) to solve online long-term constrained optimization problems where all the objectives, constraints and domains are not necessarily convex. Complexity $O(T^{8/9})$ of the algorithm is established for a proposed ESA regret which properly measures the solution quality. The application on an online extreme value-constrained river pollutant source identification, verifies the theoretical properties of our algorithm and also shows the effect of extreme value constraint in improving the solution's generalization ability. 

For future research, one possible direction is to improve the regret bound by designing another strongly concave logarithm function regularizer in the dual direction. 
Another direction is to establish 
a stronger primal direction analysis
that improves upon the $O(t\cdot T^{-1/3})$ complexity 
in Lemma~\ref{lem 6}, our results can be directly strengthened within the existing framework.
Furthermore, one can directly work with the Moreau envelope of the nonconvex objective and the Fréchet subdifferential cone of the 
feasibility
indicator. By invoking Fenchel duality, this perspective yields a substantially faster algorithm.
\bibliographystyle{unsrt}  
\bibliography{references}  

\appendix
\begin{appendices}
\section{Modified Lagrangian Oracle} \label{App A}
Assumption \ref{Assu 1} implies that both $f_{t}(x)$ and all $c_{it}(x),\,i=1,2,...,I$ are uniformly upper bounded by $f_{\max}$ and $c^i_{\max}$. Given $\epsilon >0$, define the constraint violation set $\Omega(\epsilon):=\{x\in\mathcal{X}|b_i<\frac{1}{T}\sum_{t=1}^{T} c_{it}(x)<b_i+\epsilon,\,\exists i\in\{1,2,...,I\}\}$.
\begin{pro} \label{Proposition 1}
    Consider the modified ``average'' Lagrangian first,
    \begin{equation*}
\min_{x\in\mathcal{X}}\max_{y\in[0,y_{\max}]^I}\frac{1}{T}\sum_{t=1}^{T}\left\{ f_{t}(x)+\sum_{i=1}^{I}\gamma_{i}[c_{it}(x)-b_{i}]\right\}, 
    \end{equation*}
    and $\forall\delta>0$, there exist a value $\epsilon>0$, such that $\Omega(\epsilon)$ has the Lebesgue measure smaller than $\delta$, when $y_{\max}\geq f_{\max}/\epsilon$. In addition, $\Omega(\epsilon)$ could not be $\emptyset$ if $\mathcal{X}$ is a connected set. Moreover, $\delta$ could not be $0$ if the interior of $\mathcal{X}$ 
   is still a connected set and there exits an element in the interior of $\mathcal{X}$ making at least one constraint violated.
    \label{Les}
\end{pro}

To prove Proposition \ref{Proposition 1}, we first have the following definitions for set argument. Let $V_{\epsilon}(x)$ denote the $\epsilon$-neighborhood of $x\in\mathcal{X}$ (All the distances used in this proof of are Euclidean distance which is without loss of generality). The infeasible region of the $i$-th inequality in (\ref{Lag}) is defined as $\Omega_{Ui}$ (a corresponding set is defined as $\Omega_{Ui}'$),
\[\Omega_{Ui}=\left\{x\in\mathcal{X}\Bigg|\frac{1}{T}\sum_{t=1}^T c_{it}(x)>b_i\right\},\,\Omega_{Ui}'=\left\{x\in\mathcal{X}\Bigg|\frac{1}{T}\sum_{t=1}^T c_{it}(x)<b_i\right\}.\]
Note that $\Omega_{Ui}'\neq\mathcal{X}/\Omega_{Ui}$ if $\mathcal{X}$ is a connected set. We generally assume $\Omega_{Ui}$ and $\Omega_{Ui}'$ are not empty for each $i=1,2,...,I$. And $\Omega_{Ui}$ can be parted by two disjoint sets,
\begin{align*}
    \Omega_{Ui}= & \Omega_{i}(\epsilon)\cup\Omega_{i}'(\epsilon)
    \\
    = & \left\{x\in\mathcal{X}|b+\epsilon>\frac{1}{T}\sum_{t=1}^T c_{it}(x)>b_i\right\}\cup\left\{x\in\mathcal{X}|\frac{1}{T}\sum_{t=1}^T c_{it}(x)\geq b_i+\epsilon\right\},
\end{align*}
for any $\epsilon>0$. It is known that $\Omega_{i}'(\epsilon)$ has a uniform lower bound for any $\epsilon>0$, which means this set could be paneled by a sufficiently large $y_{\max}$. For $\Omega_{Ui}$, we have the following property.
\newtheorem{prop}{Property}
\begin{prop} \label{Property 1}
For each $x\in\Omega_{Ui}$, there exists a lower bound $\epsilon_{x}$ to let \[\frac{1}{T}\sum_{t=1}^Tc_{it}(x)>b_i+\epsilon_{x}.\] When $\mathcal{X}$ is a connected set, there is no such uniform lower bound $\epsilon_u>0$ to let $\frac{1}{T}\sum_{i=1}^T c_{it}(x)>b_i+ \epsilon_u$ for all $x\in\Omega_{Ui}$. 
\end{prop}
\begin{proof}
For an arbitrary chosen $x\in\Omega_{Ui}$, we have $\frac{1}{T}\sum_{t=1}^T c_{it}(x)>b_i$. So by the density of real number, there exists a $\epsilon_x>0$ to let $\frac{1}{T}\sum_{t=1}^T c_{it}(x)>b_i+\epsilon_x>b_i$. Hence, the pointwise lower bound is proven.

By the connectivity of $\mathcal{X}$, there exists a path $\mathcal{X}^P_{x_1,\,x_2}$ from $x_1\in\Omega_{Ui}$ to $x_2\in\Omega_{Ui}'$ within $\mathcal{X}$, where function $\frac{1}{T}\sum_{t=1}^T c_{it}(x)$ is $G$-Lipschitz continuous at this path. By the Lipschitz continuous condition (Assumption \ref{Assu 1}), if there exists a uniform lower bound $\epsilon_u$ on $\Omega_{Ui}$, we have $\forall x\in \mathcal{X}^P_{x_1,\,x_2}\cap V_{\epsilon_u/G}(x^i_0)\cap\Omega_{Ui}$, such that 
\[\epsilon_u\geq G\Vert x-x^i_0\Vert\geq \frac{1}{T}\left|\sum_{t=1}^T c_{it}(x)-b_i-\sum_{t=1}^T c_{it}(x^i_0)+b_i\right|=\frac{1}{T}\sum_{t=1}^T c_{it}(x) - b_i,\]
where $x^i_0$ is a solution of $i$-th constraint function $\frac{1}{T}\sum_{t=1}^T c_{it}(x)=b_i$ to let the set $V_{\epsilon_u/G}(x^i_0)\cap\Omega_{Ui}$ is non-empty. The absolute value above is eliminated based on the definition of $\Omega_{Ui}$ above. Then by the intermediate value theorem, such a point $x^i_0$ to let $\frac{1}{T}\sum_{t=1}^T c_{it}(x)=b_i$ always exists ($\frac{1}{T}\sum_{t=1}^T c_{it}(x_2)<b_i<\frac{1}{T}\sum_{t=1}^T c_{it}(x_1)$). After denoting the set of zero points $x^i_0$ as $\mathcal{X}^P_{x_1,\,x_2,\,0}$, it generally includes two types of zero points,
\\
\\
(Type 1: Zero Point with No-sign-change)
\begin{align*}
    \biggl\{x_{0}^{i}\in\mathcal{X}_{x_{1},\,x_{2},\,0}^{P}\bigg|\exists\epsilon>0,\max_{V_{\epsilon}(x_{0}^{i})\cap\mathcal{X}_{x_{1},\,x_{2}}^{P}}\frac{1}{T}\sum_{t=1}^{T}c_{it}(x)=b_{i}\lor\min_{V_{\epsilon}(x_{0}^{i})\cap\mathcal{X}_{x_{1},\,x_{2}}^{P}}\frac{1}{T}\sum_{t=1}^{T}c_{it}(x)=b_{i}\biggr\},
\end{align*}
(Type 2: Zero Point with Sign-change)
\begin{align*}
\biggl\{ x_0^i\in\mathcal{X}^P_{x_1,\,x_2,\,0}\bigg|\forall \epsilon>0,  \max_{V_{\epsilon}(x_{0}^i)\cap\mathcal{X}^P_{x_1,\,x_2}}\frac{1}{T}\sum_{t=1}^T c_{it}(x)>b_i\land\min_{V_{\epsilon}(x_{0}^i)\cap\mathcal{X}^P_{x_1,\,x_2}}\frac{1}{T}\sum_{t=1}^T c_{it}(x)<b_i \biggr\},
\end{align*}
where Type 2 zero point $x_{0}^i$ always exists because \[\left(\frac{1}{T}\sum_{t=1}^T c_{it}(x_2)-b_i\right)\cdot\left(\frac{1}{T}\sum_{t=1}^T c_{it}(x_1)-b_i\right)<0.\]
It is obvious that points $x\in\mathcal{X}^P_{x_1,\,x_2}\cap V_{\epsilon_u/G}(x^i_0)\cap\Omega_{Ui}$ give us that $b_i + \epsilon_u\geq\frac{1}{T}\sum_{t=1}^T c_{it}(x)$. Hence, the positive uniform lower bound does not exist. 
\end{proof}

Property \ref{Property 1} illustrates that the set $\Omega_{Ui}$ has a pointwise lower bound. However, the uniform positive lower bound generally does not exist under a mild sufficient condition that $\mathcal{X}$ is a connected set. Then by $\mathcal{X}/\partial\mathcal{X}\cap V_{\epsilon_u/G}(x^i_0)\cap\Omega_{Ui}\subseteq\Omega_i(\epsilon_u)$, the set $\Omega_i(\epsilon_u)$ can not be empty for all $\epsilon_u$. Furthermore, augmented from the proof of Property \ref{Property 1}, we also suppose whether there exists a $\epsilon_u>0$ to let $\Omega_{Ui}(\epsilon_u)$ has $0$ Lebesgue measure or not. If the interior of $\mathcal{X}$, $\mathcal{X}/\partial\mathcal{X}$, is still a connected set and $\Omega_{Ui}\cap(\mathcal{X}/\partial\mathcal{X})\neq\emptyset$. Considering the path from $x_1\in\Omega_{Ui}\cap(\mathcal{X}/\partial\mathcal{X})$ to $x_2\in\Omega_{Ui}'\cap(\mathcal{X}/\partial\mathcal{X})$ within $\mathcal{X}/\partial\mathcal{X}$ and denote it as $(\mathcal{X}/\partial\mathcal{X})^P_{x_1,\,x_2}$, for $x\in (\mathcal{X}/\partial\mathcal{X})^P_{x_1,\,x_2}\cap V_{\epsilon_u/G}(x^i_0)\cap\Omega_{Ui}$, we have,
\[x\in\mathcal{X}/\partial\mathcal{X}\Rightarrow\exists \epsilon_{\partial}(x)>0,\,V_{\epsilon_{\partial}(x)}(x)\in\partial\mathcal{X}/\mathcal{X},\]
where $\epsilon_{\partial}(x)$ an index dependents on each $x\in\mathcal{X}/\partial\mathcal{X}$, so as $\epsilon_{\lambda}(x)$ and $\epsilon_{\iota}(x)$ that will be defined next.
Then considering the $x\in V_{\epsilon_u/G}(x^i_0)\cap(\mathcal{X}/\partial\mathcal{X})^P_{x_1,\,x_2}$ at the path, there exists a $\epsilon_{\lambda}(x)$ to let $V_{\epsilon_{\lambda}(x)}(x)\cap(\mathcal{X}/\partial\mathcal{X})^P_{x_1,\,x_2}\subseteq V_{\epsilon_{u}/G}(x_0^i)\cap(\mathcal{X}/\partial\mathcal{X})^P_{x_1,\,x_2}$, where $x_0^i$ is a Type 2 zero point of $(\mathcal{X}/\partial\mathcal{X})^P_{x_1,\,x_2,\,0}$, similar to the definition in the proof of Property \ref{Property 1}. By letting $\epsilon_{\iota}(x)=\min\{\epsilon_{\lambda}(x),\,\epsilon_{\partial}(x)\}$, we can let $V_{\epsilon_{\iota}(x)}(x)\subseteq V_{\epsilon_u/G}(x^i_0)\cap\Omega_{Ui}\cap(\mathcal{X}/\partial\mathcal{X})$ by the definition of interior point. If we define $m(\cdot)$ as the Lebesgue measure of a measurable set. Because $m(V_{\epsilon_{\iota}(x)}(x))>0$ and $V_{\epsilon_{\iota}(x)}\subseteq\Omega_i(\epsilon_u)$, $\Omega_i(\epsilon_u)$ could not have $0$ Lebesgue measure for all $\epsilon_u>0$. 

Then we could begin our proof of the proposition. Obviously, $m(\Omega_{Ui})=E_i<+\infty$ ($\mathcal{X}$ is a compact set). Define $\mathbb{N}^\ast$ be the set of all positive integers. We construct a series $\{\epsilon_k\}_{k\in \mathbb{N}^\ast}$ to help us formulate two set series $\{\Omega_{i}(\epsilon_k)\}_{k\in \mathbb{N}^\ast}$ and $\{\Omega_{i}'(\epsilon_k)\}_{k\in \mathbb{N}^\ast}$ to describe the extreme behavior of $\Omega_i(\epsilon)$ and $\Omega_i'(\epsilon)$ when $\epsilon\rightarrow 0$.
\begin{prop} \label{Property 2}
Given a series $\{\epsilon_k\}_{k\in \mathbb{N}^\ast}$ with three properties, (Strict positivity) $\epsilon_k>0,\,\forall k\in \mathbb{N}^\ast$,
(Strict decreasing) $\epsilon_k>\epsilon_{k+1},\,\forall k\in \mathbb{N}^\ast$,
(Convergence to zero) $\lim_{k\rightarrow +\infty}\epsilon_k=0$. Then, we could construct two set series $\{\Omega_{i}(\epsilon_k)\}_{k\in \mathbb{N}^\ast}$ with 
\[\Omega_{i}(\epsilon_k)=\left\{x\in\mathcal{X}\bigg|b_i+\epsilon_k>\frac{1}{T}\sum_{t=1}^T c_{it}(x)>b_i\right\},
\]
and $\{\Omega_{i}'(\epsilon_k)\}_{k\in \mathbb{N}^\ast}$, with
\[\Omega_{i}'(\epsilon_k)=\left\{x\in\mathcal{X}\bigg|\frac{1}{T}\sum_{t=1}^T c_{it}(x)\geq b_i+\epsilon_k\right\},\]
and with the following three properties,\\
(Containment) 
\begin{center}
$\Omega_{i}(\epsilon_{k+1})\subseteq\Omega_{Ui},\,\Omega_{i}'(\epsilon_{k})\subseteq\Omega_{Ui}$,
\end{center}
(Monotonicity)
\begin{center}
$\Omega_{i}(\epsilon_{k+1})\subseteq\Omega_{i}(\epsilon_k),\,\Omega_{i}'(\epsilon_{k})\subseteq\Omega_{i}'(\epsilon_{k+1})$,
\end{center}
(Convergence)
\begin{center}
$\cap_{j=1}^{+\infty}\cup_{k=j}^{+\infty}\Omega_{i}(\epsilon_{k})=\emptyset,\,\cup_{j=1}^{+\infty}\cap_{k=j}^{+\infty}\Omega_{i}'(\epsilon_{k})=\Omega_{Ui}$,
\end{center}
for each $i=1,2,...,I$.
\end{prop}
\begin{proof}
\emph{(Containment)} The containment is trivial by the definition of $\Omega_{i}'(\epsilon)$ and $\Omega_{i}(\epsilon)$, then for the strict positivity of each $\epsilon_k$, we finish the proof.

\emph{(Monotonicity)} By $\epsilon_k>\epsilon_{k+1},\,\forall k\in \mathbb{N}^\ast$, and considering each element $x$ in the corresponding set, we obtain,
\begin{align*}
x\in\Omega_i(\epsilon_k)
\Rightarrow  b_i+\epsilon_{k}>\frac{1}{T}\sum_{t=1}^T c_{it}(x)>b_i
\Rightarrow
b_i+\epsilon_{k+1}>\frac{1}{T}\sum_{t=1}^T c_{it}(x)>b_i
\Rightarrow  x\in\Omega_i(\epsilon_{k+1}),
\end{align*}
and
\[x\in\Omega_i'(\epsilon_{k+1})\Rightarrow \frac{1}{T}\sum_{t=1}^T c_{it}(x)>b_i+\epsilon_{k+1}\Rightarrow \frac{1}{T}\sum_{t=1}^T c_{it}(x)>b_i+\epsilon_{k}\Rightarrow x\in\Omega_i'(\epsilon_{k}).\]
However, $\Omega_{i}(\epsilon_{k+1})=\Omega_{i}(\epsilon_k),\,\Omega_{i}'(\epsilon_{k})=\Omega_{i}'(\epsilon_{k+1})$ may not hold because the set $\{x\in\mathcal{X}|b_i+\epsilon_k>\frac{1}{T}\sum_{t=1}^T c_{it}(x)>b+\epsilon_{k+1},\,k\in \mathbb{N}^\ast\}$ may not be empty.

\emph{(Convergence)} We have,
\begin{align*}
&\cap_{j=1}^{+\infty}\cup_{k=j}^{+\infty}\Omega_{i}(\epsilon_{k})=\lim_{j\rightarrow +\infty}\cup_{k=j}^{+\infty}\Omega_{i}(\epsilon_{k})=\lim_{k\rightarrow +\infty}\Omega_{i}(\epsilon_{k})
\\
=&\lim_{k\rightarrow +\infty}\left\{x\in\mathcal{X}\bigg| b+\epsilon_k>\frac{1}{T}\sum_{t=1}^T c_{it}(x)>b_i\right\},
\end{align*}
and
\begin{align*}
&\cup_{j=1}^{+\infty}\cap_{k=j}^{+\infty}\Omega_{i}'(\epsilon_{k})=\lim_{j\rightarrow +\infty}\cap_{k=j}^{+\infty}\Omega_{i}'(\epsilon_{k})=\lim_{k\rightarrow +\infty}\Omega_{i}'(\epsilon_{k})
\\
= & \lim_{k\rightarrow +\infty}\left\{x\in\mathcal{X}\bigg|\frac{1}{T}\sum_{t=1}^T c_{it}(x)\geq b_i+\epsilon_k\right\},
\end{align*}
where the second equality comes from the Monotonicity property. Suppose there exist a $x\in\Omega_{Ui}$, to let $x\in\cap_{j=1}^{+\infty}\cup_{k=j}^{+\infty}\Omega_{i}(\epsilon_{k})$ and $x\in\cup_{j=1}^{+\infty}\cap_{k=j}^{+\infty}\Omega_{i}'(\epsilon_{k})$ hold simultaneously, it leads to $\frac{1}{T}\sum_{t=1}^T c_{it}(x)>b_i$ by $x\in\Omega_{Ui}$. 
Then by Property \ref{Property 1}, there exists a $\epsilon_x$ to let, $\frac{1}{T}\sum_{t=1}^T c_{it}(x)>b_i+\epsilon_x$. And by $\lim_{k\rightarrow +\infty}\epsilon_k=0$, there exists a $K_{\epsilon_x}\in \mathbb{N}^\ast$ to let, $\forall k>K_{\epsilon_x},\,\epsilon_k<\epsilon_x$,
then $x\in\{x\in\mathcal{X}|\frac{1}{T}\sum_{t=1}^T c_{it}(x)\geq b_i+\epsilon_k\}$ and $x\notin\{x\in\mathcal{X}|b_i+\epsilon_k>\frac{1}{T}\sum_{t=1}^T c_{it}(x)> b_i\},\,\forall k>K_{\epsilon_x}$. Hence, we complete the proof of Property \ref{Property 2}.
\end{proof}

There exists series satisfies the Property \ref{Property 2}, for example, $\{1/k\}_{k\in \mathbb{N}^\ast}$. Hence, by Property \ref{Property 2} (Monotonicity), we have,
\[\lim_{j\rightarrow+\infty}m(\cup_{k=j}^{+\infty}\Omega_{i}(\epsilon_k))=m(\lim_{j\rightarrow+\infty}\cup_{k=j}^{+\infty}\Omega_{i}(\epsilon_k))=m(\cap_{j=1}^{+\infty}\cup_{k=j}^{+\infty}\Omega_{i}(\epsilon_k))=0.\]
So there exist a $j_i\in \mathbb{N}^\ast$ for each $i=1,2,...,I$ to let, 
\[m(\cup_{k=j_i}^{+\infty}\Omega_{i}(\epsilon_k))<\delta/I,\]
where $\delta$ could be any sufficiently small number larger than zero. Consequently, the complementary set,
\[\Omega_{i}'(\epsilon_{j_i})=\cap_{k={j_i}}^{+\infty}\Omega_{i}'(\epsilon_k),\]
has measure $E_i-\delta/I$. The equality holds also because by the Property 2 (Monotonicity).

All the previous arguments are about a particular $i$-th constraint. Now we try to combine them for all $i=1,2,...,I$ together. It implies that the set $\Omega_{i}'(\epsilon_{j_i})$ has a uniform lower bound $\epsilon_{j_i}$ by the definition of $\Omega_{i}'(\epsilon)$. Assigning $\epsilon=\min_{i=1,2,...,I}\{\epsilon_{j_i}\}$ and $y_{\max}=f_{\max}/\epsilon$. 
Then for all $x\in\cup_{i=1}^I\Omega_{i}'(\epsilon)$, $\frac{1}{T}\sum_{t=1}^T c_{it}(x)\geq b_i+\epsilon$.
We define $\Omega_U$ as follows,
\begin{align*}
    \Omega_U&=\left\{x\in\mathcal{X}\bigg|\frac{1}{T}\sum_{t=1}^{T} c_{it}(x)>b_i,\,i=1,2,...,I\right\}\\
    &=\cap_{i=1}^I\Omega_{Ui}=\cap_{i=1}^I(\Omega_{i}(\epsilon)\cup\Omega_{i}'(\epsilon))=(\cup_{i=1}^I\Omega_{i}(\epsilon))\cup(\cap_{i=1}^I\Omega_{i}'(\epsilon)).
\end{align*}
That's because, for all $i=1,2,...,I,$
\[x\in\Omega_{Ui}\Leftrightarrow\,\exists i=1,2,...,I,\,b_i+\epsilon_k>\frac{1}{T}\sum_{t=1}^T c_{it}(x)>b_i\,\lor\,\frac{1}{T}\sum_{t=1}^T c_{it}(x)\geq b_i+\epsilon_k.\]The above two conditions could not hold at the same time, so set $\cup_{i=1}^I\Omega_{i}(\epsilon)\cap\cap_{i=1}^I\Omega_{i}'(\epsilon)=\emptyset$. So we have, 
\begin{align*}
m(\Omega_U)= m(\cup_{i=1}^I\Omega_{i}(\epsilon))+m(\cap_{i=1}^I\Omega_{i}'(\epsilon))
 \leq  I \cdot \delta/I+m(\cap_{i=1}^I\Omega_{i}'(\epsilon))
 =\delta+m(\cap_{i=1}^I\Omega_{i}'(\epsilon)),
\end{align*}
so that $m(\cap_{i=1}^I\Omega_{i}'(\epsilon))\geq E-\delta$ (where $ m(\Omega_U)= E$) and $m(\cup_{i=1}^I\Omega_{i}(\epsilon))=\delta$ when $y_{\max}=f_{\max}/\epsilon$. If we further denote,
\[\Omega(\epsilon)=\cup_{i=1}^I\Omega_{i}(\epsilon)=\left\{x\in\mathcal{X}\bigg|b_i+\epsilon_k>\frac{1}{T}\sum_{t=1}^T c_{it}(x)>b_i,\,\exists i=1,2,...,I\right\},\]
which obviously has a Lebesgue measure smaller than $\delta$. So for all $I$ constraints, if one constraint in (\ref{Lag}) is still valid for $\mathcal{X}/\partial\mathcal{X}$, $m(\Omega_i(\epsilon))\neq 0$, then $m(\Omega(\epsilon))\neq 0$. Proposition \ref{Proposition 1} is proved.
\end{appendices}
\end{document}